\documentclass[final,onefignum,onetabnum]{siamonline220329}

\usepackage{lipsum}
\usepackage{amsfonts}
\usepackage{graphicx}
\usepackage{epstopdf}
\usepackage{algorithmic}
\ifpdf
  \DeclareGraphicsExtensions{.eps,.pdf,.png,.jpg}
\else
  \DeclareGraphicsExtensions{.eps}
\fi

\newsiamremark{remark}{Remark}
\newsiamremark{hypothesis}{Hypothesis}
\crefname{hypothesis}{Hypothesis}{Hypotheses}
\newsiamthm{claim}{Claim}

\headers{Applications of MHDM to blind
deconvolution}{T. Wolf, S. Kindermann, E. Resmerita and L. Vese}

\title{Applications of multiscale hierarchical decomposition to blind
deconvolution}

\author{Tobias Wolf\thanks{Institute of Mathematics, University of Klagenfurt, Austria 
  (\email{tobias.wolf@aau.at}, \email{elena.resmerita@aau.at}) }
\and Stefan Kindermann \thanks{Industrial Mathematics Institute, Johannes Kepler University Linz, Austria} (\email{kindermann@indmath.uni-linz.ac.at}) \and Elena Resmerita \footnotemark[1]\and  Luminita Vese \thanks{Department of Mathematics, University of California at Los Angeles (UCLA)} \email{(lvese@math.ucla.edu)}}

\usepackage{amsopn}

\makeatletter
\newcommand*{\addFileDependency}[1]{
  \typeout{(#1)}
  \@addtofilelist{#1}
  \IfFileExists{#1}{}{\typeout{No file #1.}}
}
\makeatother

\newcommand{\R}{\mathbb{R}}

\newcommand{\N}{\mathbb{N}}

\newcommand{\C}{\mathbb{C}}

\DeclareMathOperator*{\argmin}{arg\,min}
\newcommand{\abs}[1]{\left| #1 \right|}
\newcommand{\pars}[1]{\left( #1 \right)}
\newcommand{\squarebrackets}[1]{\left[ #1 \right]}

\newcommand{\seq}[1]{\pars{#1_n}_{n\in\N}}

\newcommand{\set}[1]{\left\{#1 \right\}}
\renewcommand{\Re}[1]{\text{Re}\left( #1 \right)}
\renewcommand{\Im}[1]{\text{Im}\left( #1 \right)}

\newcommand{\inv}[1]{#1^{-1}}
\newcommand{\norm}[1]{\left\Vert #1 \right\Vert }

\newcommand{\dom}{\text{dom}\,}

\newcommand{\sgn}{\text{sgn}\,}
\definecolor{unibluedark}{RGB}{68, 111, 128}
\definecolor{unibluelight}{RGB}{84, 159, 198}
\def\mathcolor#1#{\@mathcolor{#1}}
\def\@mathcolor#1#2#3{%
  \protect\leavevmode
  \begingroup
    \color#1{#2}#3%
  \endgroup
}

\usepackage{amssymb}

\usepackage{enumitem}
\usepackage{mathtools}
\usepackage[utf8]{inputenc}
\usepackage{multicol}

\usepackage{siunitx}
\usepackage[normalem]{ulem}
\allowdisplaybreaks

\ifpdf
\hypersetup{
  pdftitle={Applications of multiscale hierarchical decomposition to blind deconvolution},
  pdfauthor={Tobias Wolf, Stefan Kindermann, Elena Resmerita and Luminita Vese}
}
\fi

\begin{document}

\maketitle

\begin{abstract}
The blind image deconvolution is a challenging, highly ill-posed nonlinear inverse problem. We introduce a Multiscale Hierarchical Decomposition Method (MHDM) that is  iteratively solving variational problems with adaptive data and regularization parameters, towards obtaining finer and finer details of the unknown kernel and image. We establish convergence of the residual in the noise-free data case, and then in the noisy data case when the algorithm is stopped early by means of a discrepancy principle.  Fractional Sobolev norms are employed as regularizers for both kernel and image, with the advantage of computing the minimizers explicitly in a pointwise manner. In order to break the notorious symmetry occurring during each minimization step,  we enforce a positivity constraint on the Fourier transform of the kernels. Numerical comparisons with a single-step variational method and a non-blind MHDM show that our approach produces comparable results, while less laborious  parameter tuning is necessary at the price of more computations. Additionally, the scale decomposition of both  reconstructed kernel and image provides a meaningful interpretation of the involved iteration steps. 
\end{abstract}



\section{Introduction}

An important problem in image processing is the image restoration one, which aims to remove noise and blur from a degraded image. 
More precisely, assume that $f$ is a given blurry-noisy image, with the degradation model 
$f=k*u+n^\delta,$
where $u$ is the true image to be recovered, $k$ is a blurring kernel,
and  $n^\delta$ denotes some kind of additive noise.  There are plenty of statistical, variational and partial differential equation strategies to approach the problem. A classical variational model for this linear ill-posed problem under the assumption of normally distributed noise is
\begin{equation}\label{eq:variational}\min_{u}\{\|f-k*u\|_{L^2(\R^2)}^2+\lambda Reg(u)\},
\end{equation}
where $\lambda>0$ is the regularization parameter that should balance  stability and accuracy in the solution reconstruction, and $Reg$ stands for a penalty that promotes desired features for the recovered image, such as total variation in case of piecewise constant structures (see \cite{rudin1992nonlinear}). However, recovering both $u$ and $k$ from $f$, knowing little information about the degradation, is a highly ill-posed nonlinear inverse problem, so-called blind deconvolution. For instance, it occurs in the context of astronomical imaging \cite{Jefferies1993,Pantin2007,Prato_2013}, microscopy \cite{Kim2021,Chen2013,Avila2024} or movement correction in digital photography \cite{Cai2009,Fergus2006RemovingCS,Levin2006}. As above, one way to alleviate the difficulty in solving this problem is to use the variational approach with regularization. 

Seminal work in \cite{you1996regularization} and  \cite{chan1998total}  proposed blind 
deconvolution models using joint minimizations of the form
\begin{equation}\label{eq:conv}\min_{u,k}\{\|f-k*u\|_{L^2(\R^2)}^2+\lambda Reg(u) + \mu Reg(k)\},
\end{equation}
which can be solved using alternating minimization and two coupled Euler-Lagrange equations. Here $Reg$ denotes some generic regularization functional, which might be chosen differently for $u$ and $k$. In \cite{you1996regularization}, the regularization terms were both Sobolev $H^1$ norms, while in \cite{chan1998total}, they were the total variation for $u$, and total variation or the Sobolev norm $H^1$ for $k$. Still, the joint regularization problem, as is, involves too much symmetry between the unknowns $u$ and $k$ and thus, non-uniqueness issues appear. Note that, under some assumptions, if $(\tilde u, \tilde k)$ is a joint minimizer in \eqref{eq:conv}, then also $\left(m\tilde k,\frac{\tilde u}{m}\right)$ is a joint minimizer, for some constant $m$. A detailed analysis of this difficulty and useful ideas are presented in the book   \cite[Chapter 5]{chan2005image}. For instance, additional constraints can be included for better restoration: $k\geq0$, $\int_{\R^2} k(x)\,dx=1$, or $k$ radially symmetric, to break the symmetry of the problem. Regarding possible regularization functionals, we mention 
the works \cite{Wang2022,Shao2015,Krishnan2011} which employ sparsity promoting functionals based on $\ell^p$ (quasi)norms, as well as \cite{Mourabit2022, Shao2019} where generalizations of the total variation are used as penalty terms for the image. Note that regularization of both the image and the kernel seems necessary, since omitting a penalty on the latter might lead to inadequate results, cf \cite{Perrone2014}. Beside the variational approaches, two stage methods \cite{Polyblur,justen2006blind}, as well as a multitude of (statistically motivated)  iterative algorithms  \cite{Osher2005,Li2012,Fish95,Babacan2009,Levin2009} have been proposed. More recently, the application of machine learning methods to blind deconvolution has become popular \cite{ren2020neural,Benfenati_2023,Alban2024}.
However, in this work, we focus on a different technique, called multiscale  hierarchical decomposition of images (MHDM), introduced in \cite{VeseMultiscale, VeseDeblurring}, that favors gradual reconstruction of image features at increasingly small scales. More precisely, in that work for (denoising and) non-blind deconvolution, it was emphasized that
separating cartoon and texture in images is highly dependent on the scale
$\lambda$ from \eqref{eq:variational}, in the sense that details in an image (usually part of the texture) can be seen
as a cartoon at a refined scale, such as $\lambda/2$. 

In \cite{VeseMultiscale, VeseDeblurring}, one starts with getting a minimizer $u_0$ of \eqref{eq:variational}, then continues with iteratively solving   similar  minimization problems  which aim at extracting more detailed information from the current data $f-k*(u_0+\dots+u_{i-1)}$ by using  different scale parameters $\lambda$ at every step.  Thus, one obtains a sequence of minimizers $u_0$, $u_1$, ...,  via
$$u_i\in\arg\min_{u}\{\|f-k*(u+u_0+\dots+u_{i-1)}\|_{L^2(\R^2)}^2+\lambda_i Reg(u)\},$$
such that $f\approx k*(u_0+u_1+...+u_n)$.  Energy estimates and applications to non-blind deconvolution,  scale separation, and registration are shown in \cite{VeseDeblurring, Nachman}. Moreover, error estimates for the data-fitting term, and stopping index rules are provided in recent works \cite{MultiscaleRefinementImaging, MHDMvsTikhonov, barnett2023multiscale}, which also clearly point out that the MHDM merits are at least twofold. Namely, it provides fine recoveries of images with multiscale features that are otherwise not obtainable by single step variational models, and it is pretty robust with respect to the choice of the initial parameter $\lambda$ (and of parameters involved in the computational procedures), thus  avoiding the burden of choosing it appropriately  when performing only one step in \eqref{eq:variational}. To be fair though on the comparison, note that more computations are involved when using MHDM.

To benefit from these effects, we are interested to extend the MHDM to the more complex problem of blind deconvolution. Clearly, we do not have 
$$(k_0+k_1+...+k_n)*(u_0+u_1+...+u_n)= k_0*u_0+k_1*u_1+...+k_n*u_n,$$ 
as we would try to “blindly” apply the hierarchical decomposition method to blind deconvolution. Instead, we introduce an appropriate procedure that provides reconstructions of the kernel and of the true image of the form  in the left-hand side above. Let us first specify the notation. We consider the observed blurred and noisy image $f^\delta$ given by \begin{equation}\label{eq:image_model}
    f^\delta = K^\dagger * U^\dagger  + n^\delta,
\end{equation}
where $U^\dagger$ is the true image, $K^\dagger$ a blurring kernel, and $n^\delta$ some additive noise.  We therefore assume $U^\dagger,K^\dagger\in L^2(\R^2)$ for the convolution to be well-defined. Here $L^2(\R^2)$ denotes the space of real-valued, square integrable function on $\R^2$. The space of complex valued, square integrable functions will be denoted by$L^2(\R^2,\C)$.  Let  $J_1,J_2: L^2(\R^2) \to \R\cup\set{\infty}$ be proper, lower semicontinuous, convex and non-negative functionals. The aim  is to reconstruct $U^\dagger$ and $K^\dagger$ from the observation $f^\delta$ by adapting the MHDM to problem \eqref{eq:image_model}. That is, we would like to decompose $U^\dagger$ and $K^\dagger$ as sums\begin{equation}
    U^\dagger = \sum_{i = 0}^\infty u_i,  \quad K^\dagger = \sum_{i = 0}^\infty k_i, 
\end{equation}
where each component $k_i$ and $u_i$ contains features of $K^\dagger$ and $U^\dagger$, respectively, at a different scale. 
Let us point out that the additive decomposition of both convoluted functions has a multitude of potential applications, even beyond image deblurring. For instance, the sequence of iterates $k_i$ contains information that could be used in the classification or learning of point-spread-functions in real applications, such as in remote sensing and atmospheric sciences, where point-spread-functions are not known and have complicated features and oscillations at different scales (see for instance the GeoSTAR PSF in \cite{GeoStationaryMicrowave}). Furthermore, such decompositions could be applied to parameter identification problems with an unknown differential operator. The kernel decomposition in such problems corresponds to an additive decomposition of the problem's Green-function. Hence, the scale decomposition of this function can be interpreted as different physical effects. Those components corresponding to coarser scales describe the overall model, while the finer scales act as its refinements. 

 Therefore, we proceed as follows. Let $(\lambda_n)_{n \in \N_0}, \ (\mu_n)_{n \in \N_0}$ be decreasing sequences of positive real numbers. Additionally, let $\Phi :(L^2(\R^2)\times L^2(\R^2))\times L^2(\R^2) \to [0,\infty)$ be a measure of similarity between $K*U$ and  an image $f$, which satisfies $\Phi(K,U,f) = 0$ for all $K,U$ with $K*U = f$. 
We consider data fidelity terms for noisy observations that have the form  \begin{align}\label{eq:data_fidelities}
    &\Phi(K,U,f^\delta) = \norm{f^\delta-K*U}_{L^2}^2 +\delta_{S_1}(U) +\delta_{S_2}(K).
\end{align} Here \begin{equation*}
    \delta_S(z) = \begin{cases}0 \text{ if } z \in S\\ \infty \text{ if } z \notin S\end{cases}
\end{equation*} denotes the indicator function of a convex set $S$ and is used to encode additional assumptions such as positivity of the kernel or constraints on the means of images and kernels.

We compute the initial iterates as \begin{equation}\label{eq:MHDM_initial}
    (u_0,k_0) \in \argmin_{u,k \in L^2(\R^2)} \Phi(k,u,f^\delta) + \lambda_0 J_1(u) + \mu_0 J_2(k).
\end{equation}
Next, set $U_0 = u_0$, $K_0 = k_0$ and determine the increments $u_1,k_1$ such that $U_1 = U_0 +u_1$ and $K_1 = K_0+k_1$ via
\begin{align*}
    (u_{1},k_{1}) \in \argmin_{u,k\in L^2(\R^2)} \Bigg\{ &\norm{f^\delta-(k+K_0)*(u+U_0)}_{L^2}^2\\& +\delta_{S_1}(u+U_0) +\delta_{S_2}(k+K_0) + \lambda_{1} J_1(u) + \mu_{1} J_2(k)\Bigg\}. 
    \end{align*}
Thus, we iterate for $n \in \N_0$, \begin{equation}\label{eq:MHDM_step}
    (u_{n+1},k_{n+1}) \in \argmin_{u,k\in L^2(\R^2)} \Phi(k+K_n,u+U_n,f^\delta) + \lambda_{n+1} J_1(u) + \mu_{n+1} J_2(k),
\end{equation}
that is,
\begin{align*}
    (u_{n+1},k_{n+1}) \in \argmin_{u,k\in L^2(\R^2)}\Bigg\{& \norm{f^\delta-(k+K_n)*(u+U_n)}_{L^2}^2\\& +\delta_{S_1}(u+U_n) +\delta_{S_2}(k+K_n) + \lambda_{n+1} J_1(u) + \mu_{n+1} J_2(k)\Bigg\},   
\end{align*}
and set $U_{n+1} = u_{n+1}+U_n$, $K_{n+1} = k_{n+1}+K_n$.
Note that \eqref{eq:MHDM_step} can also be formulated as \begin{equation}\label{eq:MHDM_step_substituion}
    (U_{n+1},K_{n+1}) \in \argmin_{U,K\in L^2(\R^2)} \Phi(K,U,f^\delta) + \lambda_{n+1}J_1(U-U_n) + \mu_{n+1} J_2(K-K_n).
\end{equation}
 We extract in a nonlinear way a sequence of functions (atoms) approximating $f=K^\dagger*U^\dagger$, 
\begin{gather*}
 f\approx k_0*u_0,\\
 f\approx (k_0+k_1)*(u_0+u_1),\\
  \ldots \\ 
f\approx (k_0+k_1+...+k_n)*(u_0+u_1+...+u_n)=K_n*U_n.
\end{gather*}
This refined multiscale hierarchical blind deconvolution model will provide a better choice for the solution than the single step (variational) model, especially when reconstructing  images with different scales, as each component $(u_i,k_i)$ at a scale $(\lambda_i,\mu_i)$ contains additional information that would have been ignored at the previous, coarser scales. As usual for ill-posed problems, we stop the iterations early, according to the discrepancy principle, in order to prevent meaningless computational steps.

 We focus on Sobolev norms as regularizers for both kernel and image, and show that the iterates can be computed in a pointwise manner by means of the Fourier domain. When choosing the Fourier transforms of the kernels to be positive, we get the chance to break the unwanted symmetry which naturally occurs in such regularization frameworks. An interpretation of this choice in terms of positive definite functions is provided as well. Note that considering more modern penalties in our approach is beyond the scope of this study, since it would require accounting for a priori information and tedious computational methods specific to that setting.

 Our work is structured as follows. Section 2 presents the convergence properties of the method, as well as the stopping rule. Section 3 is dedicated to the regularization with Sobolev norms, detailing the pointwise computation of the minimizers, while Section 4 illustrates the numerical experiments that fairly compare our procedure to a single-step variational method and to a non-blind deconvolution method.

 \section{Convergence properties and stopping rule}
To the best of our knowledge, convergence results of the iterates $U_n$ and $K_n$ generated by the MHDM are not even known for simpler one-variable deblurring problems, see \cite{VeseDeblurring}.
However, convergence of the residual $f^{\delta}-K_n*U_n$ can be shown analogously to Theorems 2.1 and 2.3 in \cite{MHDMvsTikhonov}. For completeness, we will provide  a proof in the case of noisy data. 

\begin{theorem}\label{thm:Convergence_Residual}
        Let $f = K^\dagger*U^\dagger$ be the blurry, but noise-free image, and let $(U_n,K_n) = \pars{\sum_{i = 0}^n u_i,\sum_{i = 0}^n k_i}$  with $(u_i,k_i)$ attained from \eqref{eq:MHDM_initial}, \eqref{eq:MHDM_step} with $f^\delta$ replaced by $f$. Assume that $J_1$ and $J_2$ are minimal at $0$ and that there are $C_1,C_2 \ge 1$ such that $J_1(v_1-v_2) \le C_1(J_1(v_1) + J_1(v_2))$ and $J_2(h_1-h_2)\le C_2 (J_2(h_1)+J_2(h_2))$ for all $v_1,v_2,h_1,h_2$. Additionally, assume $U^\dagger \in \dom J_1$ and $K^\dagger \in \dom J_2$.  If  $\lambda_n$ and $\mu_n$ are chosen such that $\lambda_{n-1}\ge 2C\lambda_n$ and $\mu_{n-1}\ge 2C\mu_n$ for all $n\in \N$ with $C = \max\set{C_1,C_2}$, then $\Phi(U_n,K_n,f)$ is decreasing in $n$ and satisfies
        \begin{equation}\label{eq:estimate_residual}
            \Phi(U_n,K_n,f) \le \frac{2C\pars{\lambda_0 J_1(U^\dagger) + \mu_0 J_2(K^\dagger)}}{n+1}.
        \end{equation}
\end{theorem}

 It is well-known that iterative methods for ill-posed problems perturbed by noise have to be stopped early.  We  use the \textit{discrepancy principle} as a stopping criterion, namely we terminate the iteration at index $n^*(\delta)$ defined as \begin{equation}\label{eq:discrepancy_principle}
    n^*(\delta) = \max\set{n\in N : \Phi(K_n,U_n, f^\delta) >\tau \delta^2} +1
\end{equation}
for some $\tau >1$.  The well-definedness of the stopping index $n^*(\delta)$ is a consequence of the following Theorem, where for simplicity of notation the iterates obtained from the noisy observation $f^\delta$ are still denoted by $K_n$ and $U_n$.

\begin{theorem}\label{thm:noisy_residual}
    Under the assumptions of Theorem \ref{thm:Convergence_Residual}, let $U_n, K_n$ be obtained by \eqref{eq:MHDM_initial} and \eqref{eq:MHDM_step}. If  $\norm{n^\delta}_{L^2} \le \delta$ for some $\delta>0$ holds, then the residual $\Phi(K_n,U_n,f^\delta)$ is decreasing in $n$ and satisfies the estimate\begin{equation}\label{eq:estimate_residual_noisy}
         \Phi(K_n,U_n,f^\delta)\le \frac{2C\pars{\lambda_0 J_1(U^\dagger) + \mu_0 J_2(K^\dagger)}}{n+1}+ \delta^2. 
    \end{equation}
If  additionally $(n^*(\delta))_{\delta>0}$ is unbounded as $\delta\to 0$, then $(\Phi(K_{n^*(\delta)},U_{n^*(\delta)}, f^\delta))_{\delta>0}$ converges to zero.

    \begin{proof}
              Let $n_0 \in \N$. We may assume $\Phi(K_{n_0},U_{n_0},f^\delta) \ge \delta^2$, as otherwise there is nothing to show. By the optimality of $(u_{n},k_n)$ in \eqref{eq:MHDM_step}, it is \begin{equation}\label{eq:comparison_zero_noise}
                \Phi(K_{n},U_{n},f^\delta) +\lambda_{n} J_1(u_{n}) + \mu_{n}J_2(k_{n})  \le \Phi(K_{n-1},U_{n-1},f^\delta) +\lambda_{n} J_1(0) + \mu_{n}J_2(0).
            \end{equation}
        Together with by the minimality of $J_1$ and $J_2$ at $0$, this implies that $\Phi(K_n,U_n,f^\delta)$ is decreasing in $n$. On the other hand, comparing to $(U^\dagger -U_{n-1}, K^\dagger-K_{n-1})$ yields \begin{align}\label{eq:comparison_previous_noise}
                \Phi(U_{n},K_{n},f^\delta) +\lambda_{n} J_1(u_{n}) + \mu_{n}J_2(k_{n})  \notag&\le \Phi(U_{n-1} + (U^\dagger-U_{n-1}), K_{n-1} + (K^\dagger-K_{n-1}),f^\delta) \\&+  \lambda_n J_1(U^\dagger-U_{n-1}) + \mu_n J_2(K^\dagger-K_{n-1})\notag  \\&\le \delta^2 +  \lambda_n J_1(U^\dagger-U_{n-1}) + \mu_n J_2(K^\dagger-K_{n-1})
            \end{align}
            Now, for any $1 \le j \le n_0$ it is 
            \begin{equation}\label{eq:new} \Phi(K_j,U_j,f^\delta) \ge \delta^2, \end{equation} and we have \begin{align*}
           & \Phi(K_j,U_j,f^\delta) + \lambda_j J_1(U^\dagger-U_j) + \mu_j J_2(K^\dagger-K_j) \\
            &= \Phi(K_j,U_j,f^\delta) + \lambda_j J_1(U^\dagger-u_j-U_{j-1}) + \mu_j J_2(K^\dagger-k_j-K_{j-1})\\
            &\le \Phi(K_j,U_j,f^\delta) +  C \pars{\lambda_jJ_1(U^\dagger-U_{j-1}) +J_1(u_j)  + \mu_j J_2(K^\dagger-K_{j-1}) +J_2(k_j)}\\
            & = \Phi(K_{j},U_{j},f^\delta) +\lambda_j J_1(u_j) +\mu_j J_2(k_j)\\
            &+ (C-1)\pars{\lambda_j J_1(u_j) +\mu_j J_2(k_j)} \\
            & + C\pars{\lambda_j J_1(U^\dagger-U_{j-1})+\lambda_j J_2(K^\dagger-K_{j-1})}\\
            &\overset{\eqref{eq:comparison_previous_noise}}{\le} \delta^2 + (C-1)\pars{\lambda_j J_1(u_j) +\mu_j J_2(k_j)}\\
            &+(C+1)\pars{\lambda_j J_1(U^\dagger-U_{j-1})+\mu_j J_2(K^\dagger-K_{j-1})}  \\
            & \overset{\eqref{eq:comparison_previous_noise}}{\le} \delta^2 + (C-1) \pars{\delta^2 - \Phi(K_j,U_j,f) + \lambda_j J_1(U^\dagger-U_{j-1}) + \mu_j J_2(K^\dagger-K_{j-1}) } \\
            &+ (C+1) \pars{\lambda_j J_1(U^\dagger-U_{j-1})+\mu_j J_2(K^\dagger-K_{j-1})}\\
            & \overset{\eqref{eq:new}}{\le} \delta^2 + 2C\pars{\lambda_j J_1(U^\dagger-U_{j-1}) +\mu_j J_2(K^\dagger-K_{j-1})}.
        \end{align*}
                Therefore the choice $\lambda_{j-1}\ge 2C\lambda_j$ and $\mu_{j-1}\ge 2C\mu_j$ yields \begin{equation} \label{eq:help_convergence_residual_noise}
            \Phi(K_j,U_j,f^\delta) + \lambda_j J_1(U^\dagger-U_j) + \mu_j J_2(K^\dagger-K_j) \le \delta^2+ \lambda_{j-1} J_1(U^\dagger-U_{j-1}) +\mu_{j-1} J_2(K^\dagger-K_{j-1}).
        \end{equation}
            Using $0$ instead of $U_{j-1}$ and $K_{j-1}$ in the previous calculations, one obtains for $j = 0$ that  \begin{equation}\label{eq:help_convergence_residual_noisy_initial}
                \Phi(K_0,U_0,f^\delta) + \lambda_0 J_1(U^\dagger-U_0) + \mu_0 J_2(K^\dagger-k_0) \le \delta^2 + 2C(\lambda_0J_1(U^\dagger) + \mu_0 J_2(K^\dagger)).
            \end{equation}     
                Using \eqref{eq:help_convergence_residual_noise} repeatedly for $j = 1,\dots,n_0$ and employing the non-negativity of $J_1$ and $J_2$, it is \begin{align*}
            (n_0+1)&\Phi(U_{n_0}, K_{n_0}, f) \overset{\eqref{eq:comparison_zero_noise}}{\le} \sum\limits_{j = 0}^{n_0} \Phi(U_j, K_j, f) + \pars{\lambda_{n_0} J_1(U^\dagger-U_{n_0}) +\mu_{n_0} J_2(K^\dagger-K_{n_0})}  \\
           & \overset{\eqref{eq:help_convergence_residual_noise}}{\le} \delta^2+\sum\limits_{j = 0}^{n_0-1} \Phi(U_j, K_j, f) + \pars{\lambda_{n_0-1} J_1(U^\dagger-U_{n_0-1}) +\mu_{n_0-1} J_2(K^\dagger-K_{n_0-1}) }
            \\  & \le \dots\le (n_0)\delta^2 +\pars{\lambda_0 J_1(U^\dagger-U_0) + \mu_0 J_2(K^\dagger-K_0) } \\&\overset{\eqref{eq:help_convergence_residual_noisy_initial}}{\le}  (n_0+1)\delta^2 +2C\pars{\lambda_0 J_1(U^\dagger)+\mu_0 J_2(K^\dagger)}.
        \end{align*}
        Dividing by $n_0+1$ proves \eqref{eq:estimate_residual_noisy}.
    \end{proof}
\end{theorem}

\section{Sobolev norm regularizers} \label{Sobolev}
In order to illustrate our proposed method, we follow the work \cite{BurgerScherzer2001RegularizationMF} and use Sobolev norms as regularizers for kernels and images. That is, we consider the case $\Phi(K,U,f) = \norm{K*U-f}_{L^2}^2 + \delta_{S_1}(U) + \delta_{S_2}(K)$, $J_1 =  \norm{\cdot}_{H^r}^2$, and $J_2= \norm{\cdot}_{H^s}^2$ for $r,s \ge 0$. For defining the Bessel Potential norm for $r\in \R$ we  set $\Delta(x)  = 1+\abs{x}^2$. Then \begin{equation}\label{eq:definition_fractional_sobolev}
    \norm{u}_{H^r}^2 = \int\limits_{\R^2} \Delta(x)^r \abs{\hat u (x)}^2 dx
\end{equation}
where $\hat u$ is the Fourier transform of $u$.  Note that the well-definedness of the MHDM with those regularizers follows analogously to  \cite[Theorem 3.6]{justen2006blind} for all sets of constraints that are closed under addition and  satisfy $S_1 \cap \dom J_1 \neq \emptyset $ and $S_2 \cap \dom J_2 \neq \emptyset$. However, we will start with analyzing the MHDM without any constraints (i.e. $S_1 = S_2 = L^2(\R^2)$). Hence, we compute the first iterate $(U_0,K_0)$ via \begin{equation}\label{eq:H^r-H^k_initial}
    \min_{u,k\in L^2(\R^2)} \norm{k*u-f}_{L^2}^2 + \lambda_0 \norm{u}_{H^r}^2 + \mu_0\norm{k}_{H^s}^2.
\end{equation}
The norm defined in \eqref{eq:definition_fractional_sobolev} is a different, but equivalent norm than the one used in \cite{BurgerScherzer2001RegularizationMF}. Nonetheless, the results of \cite[Lemma 3.3]{BurgerScherzer2001RegularizationMF} still hold: Denote the Fourier transform of $f$ by $\hat f$ and $\bar z$ the complex conjugate of a complex number $z$. Then,  the Fourier transforms $(\hat u_0,\hat k_0)$ of all minimizers of \eqref{eq:H^r-H^k_initial} are pointwise given as
\begin{align}
    \hat u_0 = \sgn(\bar\psi \hat f) \sqrt{\squarebrackets{\sqrt{\frac{\mu_0}{\lambda_0}\Delta^{s-r}}\abs{\hat f} - \mu_0\Delta^s}_+},\label{eq:choice_u_BurSch}\\
    \hat k_0 = \sgn{(\psi)}\sqrt{\squarebrackets{\sqrt{\frac{\lambda_0}{\mu_0} \Delta^{r-s}}\abs{\hat f} - \lambda_0\Delta^r}_+}\label{eq:choice_k_BurSch},
\end{align}
for arbitrary measurable functions $\psi$ with $\psi(x) \ne 0$ for all $x$. Here $\squarebrackets{w}_+ : = \max\set{w,0}$  and \begin{equation*}
    \sgn(z) = \begin{cases}
        \frac{z}{\abs{z}}, & \text{ if } z \neq 0,\\
        0 , & \text{ if } z = 0,
    \end{cases}
\end{equation*} for $z \in \C$. From here on we assume that $f$ has limited bandwith, i.e. $\hat f$ is compactly supported. This ensures that $\hat u_0$, $\hat k_0$ as defined in \eqref{eq:choice_u_BurSch} and \eqref{eq:choice_k_BurSch} are indeed in $L^2(\R^2,\C)$ (see \cite{BurgerScherzer2001RegularizationMF}) and all subsequent minimizers can be calculated pointwise. Generally, applying the inverse Fourier transform to \eqref{eq:choice_u_BurSch} and \eqref{eq:choice_k_BurSch} gives \textit{complex-valued} minimizers of \eqref{eq:H^r-H^k_initial}. Thus, in order to obtain real-valued minimizers we need an appropriate choice for the function $\psi$. For this, recall the following result. \begin{proposition}\label{prop:characterization_realness}
    Let $u\in L^2(\R^d,\C)$ with Fourier transform $\hat u$. Then \begin{enumerate}[leftmargin = 2em]
        \item $\hat u$ is real-valued if and only if $u$ satisfies $u(x) = \overline{ u(-x)}$ for all $x \in \R^d$.
        \item $u$ is real-valued if and only if $\hat u$ satisfies $\hat u(x) = \overline{\hat u(-x)}$ for all $x \in \R^d$.
    \end{enumerate}
\begin{proof}
    See \cite[Chapter 1]{bracewell2000fourier}.
\end{proof}
\end{proposition}
Since the  image $f$ is real-valued, we obtain $\hat f(x) = \overline{\hat{f}(-x)}$ for all $x \in \R^2$. Moreover, due to the symmetry of $\Delta(x)$, we choose $\psi$ such that $\psi(x) = \overline{ \psi(-x)}$ for all $x \in \R^2$, thereby ensuring that the inverse Fourier transforms of $\hat u_0$ and $\hat k_0$ are real-valued. 
\\ For our method, we choose $\psi = 1$ and compute $u_0,\ k_0$ by applying the inverse Fourier transform. With this choice, we can interprete $\hat u_0$ as the pointwise square root of a shrinkage operator applied to $\hat f$. In particular, small frequencies in Fourier domain are eliminated, which removes noise.  Furthermore, the choice $\psi =1$ means that the Fourier transform of $k_0$ is non-negative almost everywhere. Let us also give an interpretation of this. First, recall the notion of positive semi-definite functions (see, e.g.  \cite[Section 4.4.3]{Plonka2023}).
\begin{definition}
  A function $\varphi: \R^d \to \C$ is positive semi-definite if  \begin{equation}\label{eq:psd}
  \sum_{l,m = 1}^N \varphi(x_l-x_m)\xi_l\bar \xi_m  \text{ is real-valued and satisfies }   \sum_{l,m = 1}^N \varphi(x_l-x_m)\xi_l\bar \xi_m \ge 0
\end{equation} 
for all $x_1,\dots,x_n \in \R^d$, $\xi_1,\dots,\xi_N \in \C$ and any $N \in \N$.   
\end{definition}

\begin{lemma}\label{lem:psd_and_fourier}
Let $k:\R^d \to \C$ be a function such that its Fourier transform $\hat k$ is non-negative almost everywhere. Then $k$ is positive semi-definite.\begin{proof}
    It follows analogously to the proof of \cite[Theorem 4.89]{Plonka2023}.
\end{proof}
\end{lemma}

This means, the kernel $k_0$ obtained from choosing $\psi  = 1$ in \eqref{eq:choice_k_BurSch} is positive semi-definite. In particular, it has the following properties, as one can see from \eqref{eq:psd} and Proposition \ref{prop:characterization_realness}. 
\begin{corollary}\label{cor:properties_psd} Let $\varphi: \R^d\to\C$ be positive semi-definite. Then
\begin{enumerate}[leftmargin = 2em]
    \item $\varphi(0)$ is real-valued and non-negative,      \item $\abs{\varphi(x)} \le \varphi(0)$ for all $x \in \R^d$,
    \item $\varphi(-x) = \overline{\varphi(x)} $ for all $x \in \R^d$.
\end{enumerate}
\begin{proof}
    See Chapter 12, Lemma 3 in  \cite{CheneyApproximationTheory}.
\end{proof}
\end{corollary}
Thus, $k_0$ attains a maximum at $x = 0$. Moreover, the solutions of \eqref{eq:H^r-H^k_initial} will be even functions with a peak at $x = 0$, so that they might be particularly useful to approximate conical combinations of centered Gaussians.\\
We therefore want to have the constraint   that the images have the same Fourier phase as $\hat f$ and the kernels  have non-negative Fourier-transforms for all iterates of the MHDM.  In the following we will show that the iterates of the MHDM with $S_1=S_2 = L^2(\R^2)$ can indeed be chosen to have this property.  As a result, we obtain an explicit way to compute the iterates pointwise as certain minimizers of \eqref{eq:MHDM_step} without having to enforce the constraints that $\sgn(\hat U_n) = \sgn(\hat f)$ and that $\hat K_n$ is non-negative.

\subsection*{Computation of the increments}
Let us start with  deriving how to compute the increments $(u_{n+1},k_{n+1})$ for $n\ge 0$. Generally,  we have to solve problems of the form \begin{equation}\label{eq:H^r-H^s_step}
    (u_{n+1},k_{n+1})\in \argmin_{u,k\in L^2(\R^2)}  \norm{(u+U_{n})*(k+K_{n})-f}_{L^2}^2 + \lambda_{n+1} \norm{u}_{H^r}^2 +  \mu_{n+1} \norm{k}_{H^s}^2. 
\end{equation}
Thus, in the Fourier space, this amounts to solving \begin{equation}\label{eq:H^r-H^s_step_Fourier}
\begin{split}
    (\hat u_{n+1},\hat k_{n+1})\in\argmin_{\hat u,\hat k\in L^2(\R^2,\C)} \Bigg\{ \int_{\R^2} &\abs{(\hat u(x)+\hat U_{n}(x))(\hat k(x) +\hat K_{n}(x))-\hat f(x)}^2 \\
    &+ \lambda_{n+1}\Delta(x)^r\abs{\hat u(x)}^2 +  \mu_{n+1} \Delta(x)^s\abs{\hat k(x)}^2dx \Bigg\}.
\end{split}
\end{equation}
 As for computing minimizers in the initial step \eqref{eq:H^r-H^k_initial}, we would like to solve \eqref{eq:H^r-H^s_step_Fourier} pointwise. Following the notation of \cite{BurgerScherzer2001RegularizationMF}, we fix $x\in \R^2$ and set 
\begin{align}\label{eq:substitution_H^rH^s}
    &a_{n+1} = \lambda_{n+1}\Delta^{r}(x), \ b_{n+1} = \mu_{n+1}\Delta^{s}(x),\notag \\ &q_n = \hat K_n(x),\ p_n = \hat U_n(x),\ z =\hat f(x).
\end{align} Hence, we are concerned with computing \begin{equation}\label{eq:complex_minimization_problem}
    (p_{n+1}^*,q_{n+1}^*)\in\argmin\limits_{p,q \in \C}f_n(p,q) :=\argmin_{p,q\in \C} \abs{(p+p_n)(q+q_n) -z}^2 + a_{n+1} \abs{p}^2 + b_{n+1}\abs{q}^2,
\end{equation}
 where $p_{n+1}^* = \hat u_{n+1}(x)$ and $q_{n+1}^* = \hat k_{n+1}(x)$. 
 {For the subsequent analysis of \eqref{eq:complex_minimization_problem}, we derive the first order optimality conditions with \textit{complex} variables. To this end, we consider the problem as a minimization problem in $\R^2$. We start with the following observation. For any complex number $w = \Re{w} + i\Im{w}$, we write the canonical embedding into $\R^2$ as $\tilde w = (\Re{w}, \Im{w})^T$. Hence, for $u,v,w \in \C$ it is \begin{align*}
    \abs{vu-w}_2^2 = \norm{\begin{pmatrix*}
     &\tilde v_1 &-\tilde v_2\\ &\tilde v_2 & \tilde v_1
 \end{pmatrix*}\tilde u - \tilde w}_2^2  
 \end{align*}
 Taking the gradient with respect to $\tilde u$, then gives \begin{align*}
   \nabla_{\tilde u} \pars{\norm{\begin{pmatrix*}
     &\tilde v_1 &-\tilde v_2\\ &\tilde v_2 & \tilde v_1
 \end{pmatrix*}\tilde u - \tilde w}_2^2  } =   2\begin{pmatrix*}
         &\tilde v_1 &\tilde v_2 \\ &-\tilde v_2 & \tilde v_1 \end{pmatrix*}\pars{\begin{pmatrix*}
     &\tilde v_1 &-\tilde v_2\\ &\tilde v_2 & \tilde v_1
 \end{pmatrix*}\tilde u - \tilde w} =2 \bar v (vu-w).
 \end{align*}
 Therefore, the first order optimality conditions of \eqref{eq:complex_minimization_problem} are given by \begin{align}
        (\bar q+  \bar q_n)((p+p_n)(q+q_n)-z) +a_{n+1}p = 0,\label{eq:optimality_condition_fn_1}\\
    (\bar p+\bar p_n)((p+p_n)(q+q_n)-z) + b_{n+1}q = 0.\label{eq:optimality_condition_fn_2}
\end{align}

The following Theorem shows that it is possible to choose the increment $q_n^*$ non-negative for all $n\in \N$, which therefore yields that the iterates $q_n$ are non-negative, for all $n\in \N$.
 
\begin{theorem}\label{thm:choice_q_n}
   There is a choice for $ (p_{n}^*,q_{n}^*)$ such that 
       \begin{enumerate}[leftmargin = 2em]
        \item\label{it1} $q_n \in [0,\infty)$,      
        \item\label{it2}  $\bar zp_n = z\bar p_n \in [0,\infty)$,
        \item\label{it}  $\bar z p_nq_n \le \abs{z}^2$
    \end{enumerate}
    for all $n \in \N_0$. Here $p^*_0 = p_0$ and $q^*_0 = q_0$ are understood as the solution of \eqref{eq:complex_minimization_problem} with $p_{-1} = q_{-1} = 0$.

\begin{proof}
    If $z = 0$, then $p_n = q_n = 0$ for all $n \in \N$, and the claim follows trivially. Hence, we assume $z \neq 0$. Since $a_{n+1} >0$ for all $n\in \N_0$, equation \eqref{eq:optimality_condition_fn_1} is equivalent to \begin{equation}\label{eq:critical_equation_p_general}
    p+p_n = \frac{a_{n+1}p_n +z(\bar q +  \bar q_n)}{\abs{q+q_n}^2+a_{n+1}}.  
\end{equation}
 With this, we can  prove the claims by induction.
 \begin{description}
 \item[\textit{Base case:}] Choose $q_0^* \in [0,\infty)$ according to \eqref{eq:choice_u_BurSch}. Therefore, one has $\bar q_0  = q_0 = q_0^*$ and, by \eqref{eq:critical_equation_p_general},  it follows that
 \begin{equation*}
        \bar zp_0^* = \bar zp_0 = \frac{\bar z z q_0}{\abs{q_0}^2+a_1} = \frac{q_0 \abs{z}^2}{\abs{q_0}^2 + a_1} \ge 0,
    \end{equation*}
    and \begin{equation*}
        \bar z p_0 q_0 = \frac{\abs{q_0}^2 \abs{z}^2}{\abs{q_0}^2 + a_1} = \frac{\abs{q_0}^2 }{\abs{q_0}^2 + a_1}\abs{z}^2\le \abs{ z}^2.
    \end{equation*}
    Therefore, $(i)$--$(iii)$ follow for the case $n = 0$.
   \item[\textit{Induction step:}] Assume $(i)$--$(iii)$ hold for some $n\in \N_0$. In particular, this means $\bar q_n = q_n$, so that \eqref{eq:critical_equation_p_general} becomes \begin{equation}
        \label{eq:critical_equation_p}
    p+p_n = \frac{a_{n+1}p_n +z(\bar q +  q_n)}{\abs{q+q_n}^2+a_{n+1}}.
    \end{equation} 
Hence, we can restrict the minimization of $f_n$ to minimizing over the set of all $(p,q)$, for which $p$ satisfies \eqref{eq:critical_equation_p}. We thus  need to minimize the functional \begin{align*}
    \tilde f_n(q) &:= \abs{\frac{a_{n+1}p_n(q+q_n)+z\abs{q+q_n}^2}{\abs{q+q_n}^2+a_{n+1}}-z}^2 +a_{n+1}\abs{\frac{a_{n+1}p_n +z(\bar q + q_n)}{\abs{q+q_n}^2+a_{n+1}}-p_n}^2+b_{n+1}\abs{q}^2\\&
    = \abs{\frac{a_{n+1}p_n(q+q_n) -a_{n+1}z}{\abs{q+q_n}^2+a_{n+1}}}^2 + a_{n+1} \abs{\frac{z(\bar q +q_n) - p_n\abs{q+q_n}^2}{\abs{q+q_n}^2+a_{n+1}}}^2 + b_{n+1}\abs{q}^2\\&
    = \frac{a_{n+1}^2\abs{p_n(q+q_n) -z}^2}{\pars{\abs{q+q_n}^2 +a_{n+1}}^2}  + \frac{a_{n+1}\abs{q+q_n}^2\abs{z-p_n(q+q_n)}^2 }{\pars{\abs{q+q_n}^2 +a_{n+1}}^2}+b_{n+1}\abs{q}^2\\&
    =\pars{\frac{a_{n+1}^2}{\pars{\abs{q+q_n}^2 +a_{n+1}}^2} + \frac{a_{n+1}\abs{q+q_n}^2}{\pars{\abs{q+q_n}^2 +a_{n+1}}^2}}\abs{p_n(q+q_n) -z}^2+b_{n+1}\abs{q}^2\\&
    = \frac{a_{n+1}}{\abs{q+q_n}^2 +a_{n+1}}\abs{p_n(q+q_n) -z}^2 +b_{n+1}\abs{q}^2.
\end{align*}
\begin{enumerate}[ leftmargin = 2em]
    \item In order to show that $q_{n+1} \ge 0$, we only need to show that we can choose $q_{n+1}^* \ge 0$. First, observe that $\tilde f_{n}$ must attain a minimum by coercivity and continuity. Therefore, let $\tilde q$ be a minimizer of $\tilde f_{n}$. Hence, it suffices to prove $\tilde f_n(\tilde q) \ge \tilde f_n(\abs{\tilde q})$. Let now $r\ge 0$. We show that on the set $\set{q \in \C: \abs{q} = r}$, the choice $q = r$ minimizes $\tilde f_n$. Indeed, for $q$ with $\abs{q} = r$, it holds \begin{align*}
    \tilde f_n(q) &= \frac{a_{n+1}\abs{p_n(q+q_n) -z}^2 }{\abs{q+q_n}^2 +a_{n+1}}+b_{n+1}\abs{q}^2 \\&
    = \frac{a_{n+1}\abs{p_n(q+q_n) -z}^2 }{r^2 +\abs{q_n}^2 + 2q_n \Re{q} +a_{n+1}}+b_{n+1}r^2\\&
    \ge \frac{a_{n+1}}{r^2 +\abs{q_n}^2 + 2q_n r +a_{n+1}}\abs{p_n(q+q_n) -z}^2 +b_{n+1}r^2 \\&
    = \frac{a_{n+1}\abs{p_n(q+q_n) -z}^2}{\abs{r+q_n}^2 +a_{n+1}} +b_{n+1}r^2\\&
    = \frac{a_{n+1}\abs{\bar z p_n(q+q_n)-\abs{z}^2}^2}{(\abs{r+q_n}^2 +a_{n+1})\abs{z}^2} +b_{n+1}r^2\\&
    = \frac{a_{n+1}\pars{(\bar z p_n)^2\abs{q}^2 +(\bar z p_nq_n-\abs{z}^2)^2 +2(\bar z p_n)(\bar z p_nq_n -\abs{z}^2)\Re{q}} }{(\abs{r+q_n}^2 +a_{n+1})\abs{z}^2}+b_{n+1}r^2\\&
        = \frac{a_{n+1}\pars{(\bar z p_n)^2r^2 +(\bar z p_nq_n-\abs{z}^2)^2 +2(\bar z p_n)(\bar z p_nq_n -\abs{z}^2)\Re{q}} }{(\abs{r+q_n}^2 +a_{n+1})\abs{z}^2}+b_{n+1}r^2\\&
    \ge \frac{a_{n+1}\pars{(\bar z p_n)^2r^2 +(\bar z p_nq_n-\abs{z}^2)^2 +2(\bar z p_n)(\bar z p_nq_n -\abs{z}^2)r} }{(\abs{r+q_n}^2 +a_{n+1})\abs{z}^2}+b_{n+1}r^2\\&
    = \frac{a_{n+1}\abs{\bar z p_n(r+q_n)-\abs{z}^2}^2}{(\abs{r+q_n}^2 +a_{n+1})\abs{z}^2} +b_{n+1}r^2 \\&= \tilde f_n(r),
\end{align*}
where the last inequality follows from $\bar z p_n \ge 0$ and $\bar z p_nq_n \le \abs{z}^2$.
Thus, for any $r\ge 0$ the minimum of $\tilde f_n$ on the circle of radius $r$ is attained for $q = r$. In particular, we obtain $\tilde f_n(\tilde q) \ge \tilde f_n(\abs{\tilde q})$.
\item From the previous part, we know that there is a pair $(p_{n+1}^*,q_{n+1}^*)$ minimizing $f_n$ such that $q_{n+1}^* \ge 0$. By definition, it is $p_{n+1} = p_{n+1}^*+p_n$, so we can multiply \eqref{eq:critical_equation_p} with $\bar z$ to obtain \begin{align*}
    \bar z p_{n+1} = \bar z \frac{a_{n+1}p_n+z (q_{n}^*+q_n)}{\abs{q_n^*+q_n}^2 +a_{n+1}} = \frac{a_{n+1} \bar z p_n + \abs{z}^2 q_{n+1}}{\abs{q_{n+1}}^2 +a_{n+1}} \ge 0.
\end{align*}  
\item If $p_{n+1} = 0$, the claim follows trivially. Thus, we assume $p_{n+1} \neq 0$ and find from \eqref{eq:optimality_condition_fn_2} that \begin{equation*}
    \bar p_{n+1}(q_{n+1}p_{n+1}-z) + b_{n+1}q_{n+1}^* = 0.
\end{equation*}
 Multiplying with $\bar z$ and rearranging implies \begin{equation*}
    \bar z p_{n+1} q_{n+1} = -b_{n+1}q_{n+1}^* \frac{\bar z}{\bar p_{n+1}} +\abs{z}^2.
\end{equation*}
Denote $\varphi_z = \arg z$ and $\varphi_{p_{n+1}} = \arg{p_{n+1}}$. It is \begin{equation*}
    \frac{\bar z }{\bar p_{n+1}} = \abs{z}\inv{\abs{p_{n+1}}} e^{-i\varphi_z} \inv{\pars{e^{-i\varphi_{p_{n+1}}}}} = \frac{\abs{z}}{\abs{p_{n+1}}} e^{i(\varphi_{p_{n+1}}-\varphi_z)} \ge 0,
\end{equation*}
since $\varphi_{p_{n+1}} - \varphi_z = \arg \pars{\bar z p_{n+1}}$. Therefore, $ -b_{n+1}q_{n+1}^* \frac{\bar z}{\bar p_{n+1}} \le 0$, which yields $\bar z p_{n+1} q_{n+1} \le \abs{z}^2$.
\end{enumerate}

\end{description}
\end{proof}
\end{theorem}

 {
Note that we can use $q_{n+1}^*$ instead of $q$ in \eqref{eq:critical_equation_p_general} to compute $p_{n+1}^*$. Consequently $\hat U_{n+1}$ is obtained pointwise via
\begin{equation}\label{eq:critical_equation_U_general}
    \Hat U_{n+1} = \frac{(\lambda_{n+1}\Delta^s) \hat U_n +\hat f \hat K_{n+1}}{\abs{\hat K_{n+1}}^2 + \lambda_{n+1}\Delta^s}.  
\end{equation}
However, this does not guarantee that applying the inverse Fourier transform to $\hat U_{n+1}$ and $\hat K_{n+1}$ yields real-valued solutions of \eqref{eq:H^r-H^s_step}. In order to find such solutions we make use of the following observation.
{\begin{lemma}\label{lemma:real_valued_step}
    Assume that $U_n$ and $K_n$ are real valued and let $x \in \R^2$. If $(u^*,k^*)$ is a solution of \begin{align*}
        (u^*,k^*) \in \min\limits_{u,k  \in \C}\bigg\{&\abs{\hat f(x) - (u+\hat U_n(x))(k + \hat K_n(x))}^2 \\&+ \lambda_{n+1}\Delta(x)^r\abs{u}^2 + \mu_{n+1}\Delta(x)^s\abs{k}^2\bigg\},
    \end{align*}
    then $(\bar u^*,\bar k^*)$ is a solution of \begin{align*}
        (u^*,k^*) \in \min\limits_{u,k \in \C}\bigg\{&\abs{\hat f(-x) - (u+\hat U_n(-x))(k + \hat K_n(-x))}^2 \\&+ \lambda_{n+1}\Delta(-x)^r\abs{u}^2 + \mu_{n+1}\Delta(-x)^s\abs{k}^2\bigg\}.
    \end{align*}
    \begin{proof}
    Note that $\hat f(x) = \overline{\hat f(-x)}$, $\hat U_n(x) = {\overline{{\hat{U}}_n(-x)}}$, and $\hat K_n(x)  = {\overline{{\hat{K}}_n(-x)} }$ by Proposition \ref{prop:characterization_realness}. Furthermore, we have $\Delta(x)  = \Delta(-x)$.
        Let  \begin{equation*}
            h(u,k) := \abs{\hat f(x) - (u+\hat U_n(x))(k + \hat K_n(x))}^2 + \lambda_{n+1}\Delta(x)^r\abs{u}^2 + \mu_{n+1}\Delta(x)^s\abs{k}^2,
        \end{equation*}
        and \begin{equation*}
            \tilde h(u,k) := \abs{\hat f(-x) - (u+\hat U_n(-x))(k + \hat K_n(-x))}^2 + \lambda_{n+1}\Delta(-x)^r\abs{u}^2 + \mu_{n+1}\Delta(-x)^s\abs{k}^2.
        \end{equation*}
        Since it holds that \begin{align*}
            \tilde h(\bar u,\bar v) &= \abs{\hat f(-x) - (\bar u+\hat U_n(-x))(\bar k + \hat K_n(-x))}^2 + \lambda_{n+1}\Delta(-x)^r\abs{\bar u}^2 + \mu_{n+1}\Delta(-x)^s\abs{\bar k}^2 \\ &= \abs{\overline{\hat f}(x) - (\bar u+{\overline{\hat{U}_n(x)}})(\bar k + \overline{\hat{K}_n(x)})}^2 + \lambda_{n+1}\Delta(x)^r\abs{ u}^2 + \mu_{n+1}\Delta(x)^s\abs{k}^2\\&= \abs{\hat f(x) - (u+\hat U_n(x))(k + \hat K_n(x))}^2 + \lambda_{n+1}\Delta(x)^r\abs{u}^2 + \mu_{n+1}\Delta(x)^s\abs{k}^2 = h(u,k),
        \end{align*}
        the claim follows.
    \end{proof}
\end{lemma}
Combining Theorem \ref{thm:choice_q_n} and Lemma \ref{lemma:real_valued_step} shows that we can choose the sequence $(U_n,K_n)_{n \in \N_0}$ such that $U_n$ and $K_n$ are real valued, and $\hat K_n$ is non-negative for every $n \in \N$. We summarize this in the following Lemma. 
{\begin{lemma}\label{lemma:existene_real_solutions}
    There are sequences $(U_n,K_n)_{n \in \N_0}$ of real-valued functions with
   \begin{align*}
       (U_0,K_0) \in \argmin\limits_{u,k \in L^2(\R^2)}\set{\norm{f-k*u}_{L^2}^2 + \lambda_0 \norm{u}^2_{H^r} + \mu_0 \norm{k}_{H^s}^2}
   \end{align*} 
   and 
   \begin{align*}
       (U_{n+1}-U_{n} , K_{n+1}-K_{n}) \in \argmin\limits_{u,k \in L^2(\R^2)} \bigg\{&\norm{f-(k+K_n)*(u+U_n)}_{L^2}^2 \\&+ \lambda_{n+1} \norm{u}^2_{H^r} + \mu_{n+1} \norm{k}_{H^s}^2\bigg\}
   \end{align*}
   for $n \in \N_0$, such that $\hat K_n$ is real-valued and non-negative. Moreover $(U_n,K_n)$ can iteratively be computed as follows:\begin{enumerate}[leftmargin = 2em]
       \item    $U_0$ and $K_0$ are the inverse Fourier transforms of $\hat U_0$ and $\hat K_0$ in \eqref{eq:choice_u_BurSch} and \eqref{eq:choice_k_BurSch}, respectively, with $\Psi = 1$.
       \item For $n \in \N_0$, $U_{n+1}$ and $K_{n+1}$ are the inverse Fourier transforms of $\hat U_{n+1} = \hat U_n + \hat u_{n+1}$ and $\hat K_{n+1} = \hat K_n +\hat k_{n+1}$, where $\hat u_{n+1}$, $\hat k_{n+1}$ solve \begin{align*}
           (\hat u_{n+1},\hat k_{n+1}) \in \argmin\limits_{\hat u,\hat k \in L^2(\R^2,\C)}\bigg\{&\norm{\hat f - (\hat u+\hat U_n)(\hat k+\hat K_n)}_{L^2}^2 \\&+ \lambda_{n+1} \norm{\Delta^r \hat u}_{L^2}^2 + \mu_{n+1} \norm{\Delta^s\hat k}_{L^2}^2\bigg\}
       \end{align*}
       and satisfy $\hat k_{n+1}(x) \in [0,\infty)$, $\hat k_{n+1}(x) = \hat k_{n+1}(-x)$, $\hat u_{n+1}(x) = \overline{\hat{u}_{n+1}(-x)}$ for all $x \in \R^2$, as well as $\sgn(\hat u_{n+1}) = \sgn(\hat f)$.
   \end{enumerate}
   \begin{proof}  
       From Lemma \ref{lemma:real_valued_step}, we know that it suffices to find  sequences $(\hat U_n,\hat K_n)_{n \in \N}$ of pointwise minimizers of the integrand in \eqref{eq:H^r-H^s_step} on the half-plane $\in [0,\infty) \times \R$. The existence of such $\hat K_n$ follows from Theorem \ref{thm:choice_q_n}. In particular, we have $\hat K_n(x) \in [0,\infty)$ so that by induction over \eqref{eq:critical_equation_U_general} we obtain $\hat U_n$ with $\sgn (\hat U_n(x)) = \sgn(\hat f(x)) $. Now, extend $\hat U_n$ and $\hat K_n$ by setting $\hat K_n(x) := \hat K_n(-x)$ and $\hat U_n(x) := \overline{\hat{U}_{n}(-x)}$ on $(-\infty,0]\times \R$. According to Lemma \ref{lemma:real_valued_step} these extended functions are minimizers of \eqref{eq:H^r-H^s_step_Fourier}. The realness of the inverse Fourier transforms is a consequence of Poposition \ref{prop:characterization_realness} and the fact that $f$ is real-valued.
   \end{proof}
\end{lemma}

  We can now construct a method to directly solve \eqref{eq:complex_minimization_problem}. Following the proof of Lemma \ref{lemma:existene_real_solutions}, it suffices to solve \eqref{eq:complex_minimization_problem} pointwise for $x \in [0,\infty)\times \R$. For better readability, we again use the notation \eqref{eq:substitution_H^rH^s}. First, we plug \eqref{eq:critical_equation_p} in \eqref{eq:optimality_condition_fn_2} to obtain  \begin{align*}
    \abs{\frac{a_{n+1}p_n +z(\bar q +q_n)}{\abs{q+q_n}^2+a_{n+1}} }^2(q+q_n) -z\frac{a_{n+1}\bar p_n +\bar z( q + q_n)}{\abs{q+q_n}^2+a_{n+1}}  +b_{n+1}q = 0.
\end{align*}Multiplying with $(\abs{q+q_n}^2+a_{n+1})^2$, we get
\begin{align}\label{eq:q+q_n}
    &\abs{a_{n+1}\bar p_n +\bar z (q+q_n)}^2(q+q_n)\notag \\
    &\qquad \qquad - (a_nz\bar p_n+\abs{z}^2(q+q_n))(\abs{q+q_n}^2+a_{n+1})\notag \\
    &\qquad \qquad  + b_{n+1}q(\abs{q+q_n}^2+a_{n+1})^2 = 0.
\end{align}
Since we are interested in finding a real solution of this equation, we can assume $q \in \R$ and expand to \begin{align}\label{eq:q+q_n_expanded}
    &b_{n+1}(q+q_n)^5 -b_{n+1}q_n(q+q_n)^4 + 2a_nb_{n+1}(q+q_n)^3 \notag\\&
    \qquad \qquad + \squarebrackets{2a_{n+1} \text{Re}(\bar p_n z) -a_nz\bar p_n-2a_nb_{n+1}q_n}(q+q_n)^2\notag\\
    &\qquad \qquad +\squarebrackets{a_{n+1}^2\abs{p_n}^2-a_{n+1}\abs{z}^2+a_{n+1}^2b_{n+1}}(q+q_n)-a_{n+1}^2\pars{z\bar p_n +b_{n+1}q_n}=0.
\end{align}
Since both $q_n$ and $z \bar p_n$ are real numbers, this is a polynomial of degree $5$ with only real coefficients. We can thus choose $q_{n+1}^*$ to be the real root of \eqref{eq:q+q_n_expanded}, for which $f_n$ as in \eqref{eq:complex_minimization_problem} attains the smallest value.  Then  $p_{n+1}^*$ is determined by \eqref{eq:critical_equation_p_general}. We have therefore computed $\hat u_{n+1}$ and $\hat k_{n+1}$ on $[0,\infty)\times \R$. Extending with complex conjugation to $\R^2$ and taking the inverse Fourier transform, then gives real-valued solutions of \eqref{eq:H^r-H^s_step}, such that the $\hat k_{n+1}$ is real-valued and non-negative. Thus, by choosing the those minimizers for which the Fourier transform of the kernel is real and non-negative in the MHDM iteration,} we  implicitly incorporate the constraint $\hat K_n(x) \ge 0$ for all $x$. Therefore, all approximations of the kernel will be positive semi-definite enforcing properties $(i)$--$(iii)$ in Corollary \ref{cor:properties_psd}, that acts as a way to break the symmetry of problem \eqref{eq:H^r-H^k_initial}.

\begin{remark}
    If more information on the phase of the true kernel is available, one might want to  one choose $\psi$ in \eqref{eq:choice_u_BurSch} and \eqref{eq:choice_k_BurSch} to be a different complex-valued function instead of a non-negative one. In this case, the sequence $(\hat u_n,\hat k_n)_{n \in \N}$ generated by \eqref{eq:H^r-H^s_step_Fourier} can be chosen such that $\sgn{\hat k_n} = \sgn{\psi}$ for all $n\in \N$. To see this,
    notice that substituting $\sgn{(\psi)}\hat k$, $\sgn{(\psi)}\hat K_n$, $\sgn{(\bar \psi\hat f)}\hat u$ and $\sgn{(\bar \psi\hat f)}\hat U_n$ for $\hat k$, $\hat K_n$, $\hat u$ and $\hat U_n$, respectively, in \eqref{eq:H^r-H^s_step_Fourier}, does not change the value of the objective function. Thus, we can use the iterates obtained with the constraint $\hat k \ge 0$ to compute the iterates of a MHDM with the constraint $\sgn(\hat k) = \sgn{(\psi)}$.
\end{remark}

In order to further ensure that the iterates are reasonable approximations of the true image and kernel, we additionally assume $f\in L^1$ and impose constraints on their means: 
\begin{equation}\label{eq:mean_constraints}
\int_{\R^2} K_n = 1, \quad \int_{\R^2}U_n = \int_{\R^2}f
\end{equation}

 for all $n\in \N$. Those constraints are fairly standard and can for instance be found in \cite{chan2005image, Osher2005}. In summary, this means we are considering the constraint sets:\begin{equation*}
    S_1 = \set{u:\R^2 \to \R : \int_{\R^2} u  = \int_{\R^2} f}, \quad S_2 = \set{k: \R^2 \to \R: \int_{\R^2} k =  1,\  \hat k \ge 0}.
\end{equation*}
It is not clear if these constraints can be translated into Fourier space such that the iterates of the MHDM can still be computed in a pointwise manner. Since for any function $u \in L^1(\R^2)$   it is known that $\hat u$ is uniformly continuous \cite[p.1]{bochner1949fourier}, it must be $\int_{\R^2} u = \hat u(0)$. However, simply enforcing the constraint $\hat U_{n} (0) = \int_{\R^2} f$ in \eqref{eq:H^r-H^s_step_Fourier} does not affect $\int_{\R^2} U_{n}$ as the resulting minimizer in Fourier space would only differ from the unconstrained one on a set of measure $0$ and $\hat U_n$ is an element of $H^r(\R^n)$ and not necessarily of $L^1(\R^n)$. Note that this problem does not occur in the discretized setting we use for the numerical experiments, as will be outlined later.

\begin{remark}
   Instead of employing a Bessel Potential norm as a penalty term for the image, it would be a natural idea to use a functional that favors expected structures. For instance, one could use the total variation to promote cartoon-like images. However, numerical experiments suggest that the iterates of such blind-deconvolution MHDM seem to approximate the trivial solution $u =f$ and $k =\delta$, where $\delta$ denotes the Dirac delta distribution. This could possibly be overcome by a specific choice for the sequences $\seq{\lambda}$ and $\seq{\mu}$, but would first require a deeper analysis of convergence behavior of the blind MHDM in this scenario, which is not within the scope of this work.  
\end{remark}

\section[Numerical Experiments]{Numerical Experiments\footnote{The program code is available as ancillary file from the arXiv page of this paper.}}\label{section:numerics}

The goal of our numerical experiments is to illustrate the behavior and robustness of the blind deconvolution MHDM. To this end, we compare the reconstructed image and kernel from the proposed method to those obtained from using a non-blind MHDM or a single step variational regularization, as in \eqref{eq:H^r-H^k_initial}. To achieve a fair comparison, we use the same regularizers for all methods under investigation. By using squared Bessel Potential norms, we obtain reasonable reconstructions that might not necessarily outperform methods with more problem specific regularizers. However, we think that comparing our approach to such methods should be done with a version of the MHDM that also uses more sophisticated regularizers, that is beyond the scope of this work. We want to stress that the advantage of our method is the interpretability of the scale decomposition and the potential to adapt it to a multitude of regularizers that possibly could vary along the iterations rather than significantly better reconstructions than those of a single-step variational approach with optimal regularization parameters.

\subsection{Discretization and Implementation}
Recall that  $J_1 = \norm{\cdot}_{H^r}^2$ and $J_2 = \norm{\cdot}_{H^s}^2$ for $r,s\ge 0$. We discretize an image $u$ supported on a rectangular domain $ (a,b)\times(c,d)$ as a matrix, i.e., $u \in \R^{m\times n}$. Following the derivation in Section 7.1.2 of \cite{justen2006blind}, we define the weight matrix for the Sobolev norm in Fourier space by the matrix  $\Delta \in \R^{m\times n}$ with entries \begin{equation*}
    \Delta_{i,j} = 1+2m^2\pars{1-\cos\pars{\frac{2\pi i}{m}}} + 2n^2\pars{1-\cos\pars{\frac{2\pi j}{n}}}.
\end{equation*}
Hence, a discretization of the Sobolev norm is given by \begin{equation*}
    \norm{u}_{H^r}^2 = \sum\limits_{i = 1}^m\sum\limits_{j = 1}^n \Delta_{i,j}^r \abs{\hat u_{i,j}}^2,
\end{equation*}
where $\hat u$ denotes the discrete $2D$ Fourier transform. Therefore, the $n$-th step of the MHDM is given by the pointwise update rule \begin{equation}\label{eq:MHDM_step_H^rH^s_discrete}
(\hat k_{i,j}^{(n)},\hat u_{i,j}^{(n)}) \in \argmin_{\hat k,\hat u \in \C} \set{\pars{(\hat k + \hat k_{i,j}^{(n-1)})(\hat u+\hat u_{i,j}^{(n-1)}) -\hat f_{i,j}}^2 + \lambda_n \Delta_{i,j}^r\abs{\hat u}^2  +   \mu_n  \Delta_{i,j}^s\abs{\hat k}^2},
\end{equation}
 with $\hat k_{i,j}^{-1} = u_{i,j}^{-1} = 0$. Thus for $n = 0$, we can solve \eqref{eq:MHDM_step_H^rH^s_discrete} by 
    \begin{align*}
    &\hat u_{i,j}^{(0)} = \sgn( \hat f_{i,j}) \sqrt{\squarebrackets{\sqrt{\frac{\mu_0}{\lambda_0}\Delta^{s-r}_{i,j}}\abs{\hat f_{i,j}} - \mu_0\Delta^s_{i,j}}_+}\\
    &\hat k_{i,j}^{(0)} = \sqrt{\squarebrackets{\sqrt{\frac{\lambda_0}{\mu_0} \Delta_{i,j}^{r-s}}\abs{\hat f_{i,j}} - \lambda_0\Delta^r_{i,j}}_+}.
\end{align*}
For $n\ge 1$, we make the  same substitutions as in \eqref{eq:substitution_H^rH^s}:
\begin{align*}
        &a_n= \lambda_{n}\Delta^{r}_{i,j}, \ b_n = \mu_{n}\Delta^{s}_{i,j},\notag \\ &q_n = \hat k_{i,j}^{(n-1)},\ p_n = \hat u_{i,j}^{(n-1)},\ z =\hat f_{i,j},
\end{align*}
To find critical pairs $(\hat k_{i,j},\hat u_{i,j})$, we compute the positive roots of \eqref{eq:q+q_n_expanded} that yield candidates for $\hat k_{i,j}$ and use \eqref{eq:critical_equation_p} to obtain the corresponding candidates for $\hat u_{i,j}$. The minimizing pair can then be found by choosing the critical pair, which gives the smallest value for the objective function in \eqref{eq:MHDM_step_H^rH^s_discrete} and enforce the antiysmmetry condition as described in the proof of Lemma \ref{lemma:existene_real_solutions}. In order to obtain meaningful reconstructions $u^{(n)},k^{(n)}$, we employ the additional constraint \eqref{eq:mean_constraints} that the mean of the kernel is $1$ and the mean of the reconstructed image matches the mean of the observation. In the discretization that means  \begin{equation*}
    \sum\limits_{i = 1}^m\sum\limits_{j=1}^n k_{i,j}^{(n)} = 1,\quad \sum\limits_{i = 1}^m\sum\limits_{j=1}^n u_{i,j}^{(n)} = \sum\limits_{i = 1}^m\sum\limits_{j=1}^n f_{i,j}^{(n)}.
\end{equation*}
Since we are using discrete Fourier transforms, these constraints are equivalent to $\hat k_{i,1} = 1$ and $\hat u_{1,1} = \sum\limits_{i = 1}^m\sum\limits_{j=1}^n f_{i,j}^{(n)}$. Thus, we implement them by making the updates \begin{equation*}
    \hat k_{1,1}^{(n)} = \begin{cases}
        1 &\text{ if } n = 0,\\
        0 &\text{ if } n \ge 1.
    \end{cases}\quad\text{ and }\quad  
        \hat u_{1,1}^{(n)} = \begin{cases}
            \hat f_{1,1} &\text{ if } n = 0,\\
            0 &\text{ if } n \ge 1.
        \end{cases}
\end{equation*}
instead of the previous procedure for the first entries of the matrices.

\subsection{ Behavior of the blind-deconvolution MHDM}
 For illustrating the behavior of the proposed method, we consider the image \textit{Barbara} (denoted by $U^\dagger$), which has been blurred by convolution with two different kernels. For the first blurring we choose a Gaussian kernel $K_1^\dagger$ of mean $\mu =0$ and variance $\sigma = 8$. In the second kernel $K_2^\dagger$, we use a convex combination of several Gaussians.
In both cases, the blurred image was additionally corrupted with additive Gaussian noise (mean $\mu = 0$, variance $\sigma = 4\times10^{-4}$). That is, we deal with observations $f_1^\delta,\ f_2^\delta$ obtained via $f_i^\delta = U_\text{true}*K_i + n^\delta$ for $i =1,2$. 
The true image, blurred images and noise corrupted blurred images can be found in Figure \ref{fig:H^r-H^s-input}, the corresponding kernels are shown in Figure \ref{fig:H^r-H^s-kernels}.

\begin{figure}
\centering
\begin{minipage}[c]{0.3\textwidth}
    \includegraphics[scale = 0.5]{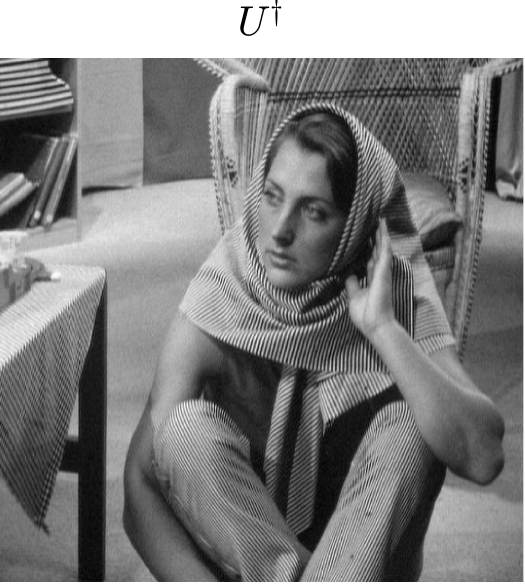}
\end{minipage}
\begin{minipage}[c]{0.3\textwidth}
 \includegraphics[scale = 0.5]{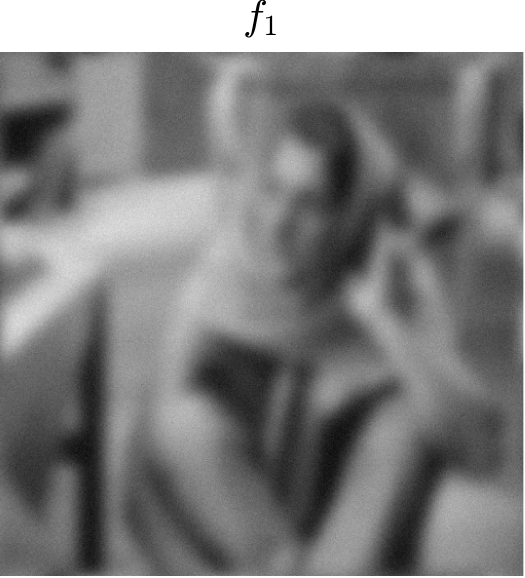}
\end{minipage}
\begin{minipage}[c]{0.3\textwidth}
    \includegraphics[scale = 0.5]{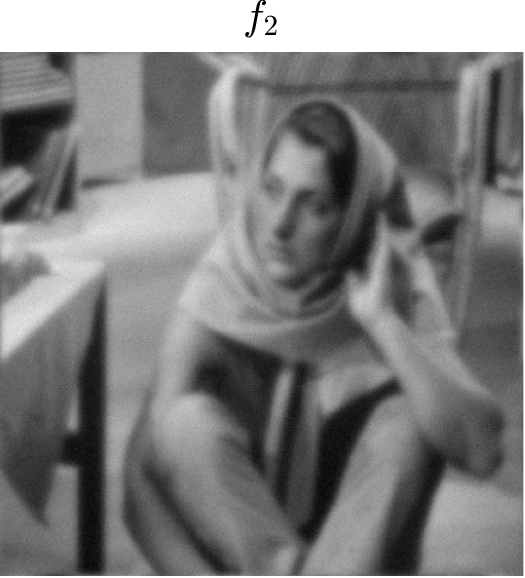}
\end{minipage}
    \caption{From left to right: true test image $U^\dagger$, observation obtained with single Gaussian blur and noise $f_1^\delta$ and observation obtained with mixture of Gaussians and noise $f_2^\delta$.}
    \label{fig:H^r-H^s-input}
\end{figure}

\begin{figure}
    \centering
\begin{minipage}[c]{0.49\textwidth}
  \includegraphics[scale = 0.49]{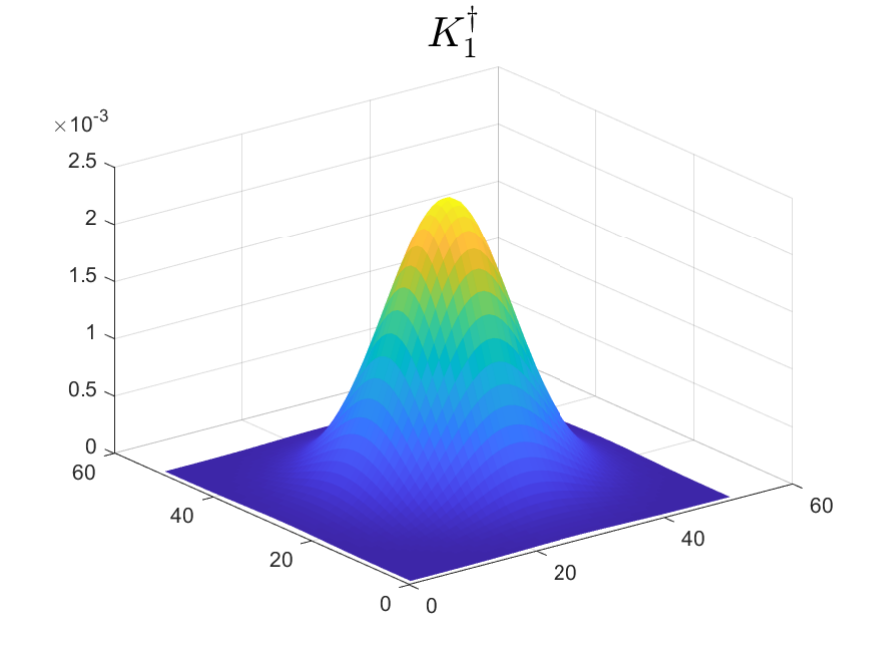}
\end{minipage}
\begin{minipage}[c]{0.5\textwidth}
 \includegraphics[scale = 0.5]{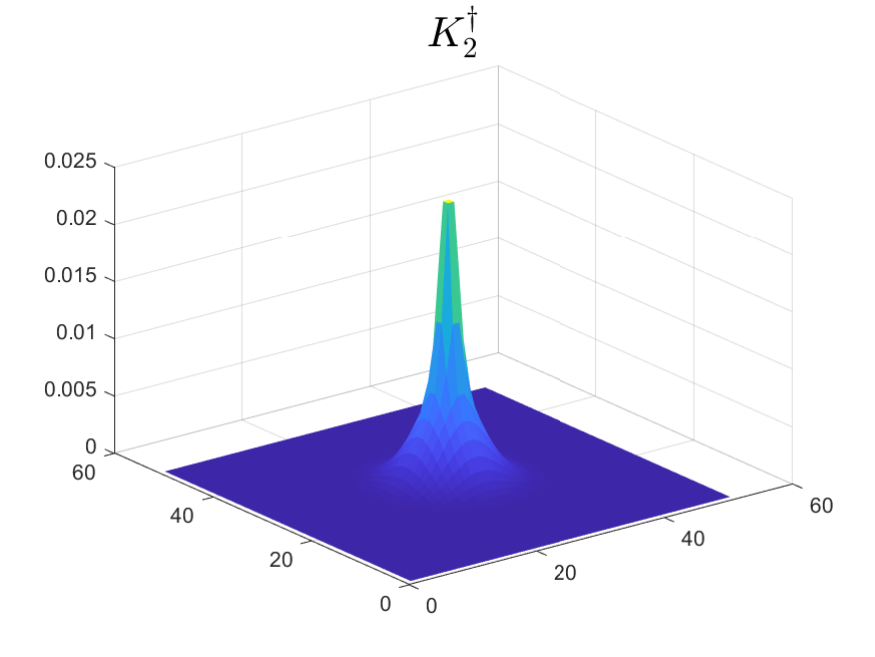}
\end{minipage}
    \caption{Left: single Gaussian kernel used to obtain $f_1^\delta$, Right: Convex combination of multiple Gaussian kernels used to obtain $f_2^\delta$.}
    \label{fig:H^r-H^s-kernels}
\end{figure}

 Curiously, the results of all our experiments (including the more extensive experiments later) improve if $s$ is chosen smaller than $r$. This means, we penalize the image with a more smoothness promoting regularizer than the kernel. Furthermore, we  observe that as long as the ratio $\frac{\lambda_0}{\mu_0}$ is constant, the actual choice of the initial parameters $\lambda_0,\mu_0$ does not significantly influence the quality of the final iterates but only the number of iterations needed until the discrepancy principle is satisfied.  \\  
 For our experiments we choose $r = 1$ and $s = 0.1$, and run the MHDM with initial parameters $\lambda_0 = 1.4\times 10^{-4}$, $\mu_0 =6.3\times 10^5 $. In accordance with Theorem \ref{thm:noisy_residual}, we choose the parameters at the $n$-th step as $\lambda_n = 4^{-n}\lambda_0$ and $\mu_n = 4^{-n} \mu_0$. Since in this experiment we artificially added noise and hence know the exact noise level, the iteration is stopped according to the discrepancy principle \eqref{eq:discrepancy_principle} with $\tau = \sqrt{1.001}$.
 
The reconstructed kernels and images are illustrated in Figure \ref{fig:blind_reconstructions}. Figure \ref{fig:H^r-H^s-blind_residual} shows the decay of the residual for both experiments. In both figures, one can clearly see the monotone decrease of the residual, confirming the theoretical results from Theorem \ref{thm:Convergence_Residual}. Figure \ref{fig:H^r-H^s-scales_image} shows the different scales $u_n$ that are obtained with the MHDM employed for the  observation $f_2^\delta$. One can see that each step adds another layer of details to the reconstruction. The corresponding scales $k_n$ and iterates $K_n = \sum\limits_{i = 0}^n k_n$  for the reconstructed kernel can be seen in Figures \ref{fig:H^r-H^s-scales_kernel} and \ref{fig:H^r-H^s-iterates_kernel}. It appears that the role of the scales is to adapt the reconstructed kernel in a twofold way. On the one side,  the height of the peak seems to increase along the iterations, while its radius decreases. On the other side, the  off-peak oscillations appear to be flattend. Notably, the early coarse scales seem to recover the  support of the bump, and the fine scales mostly shape the height of it. In the experiment with data $f_1^\delta$, the scale decomposition for the reconstructed image and kernel look similar.

\begin{figure}
    \centering
    \begin{minipage}[c]{0.49\textwidth}
   \includegraphics[scale = 0.4]{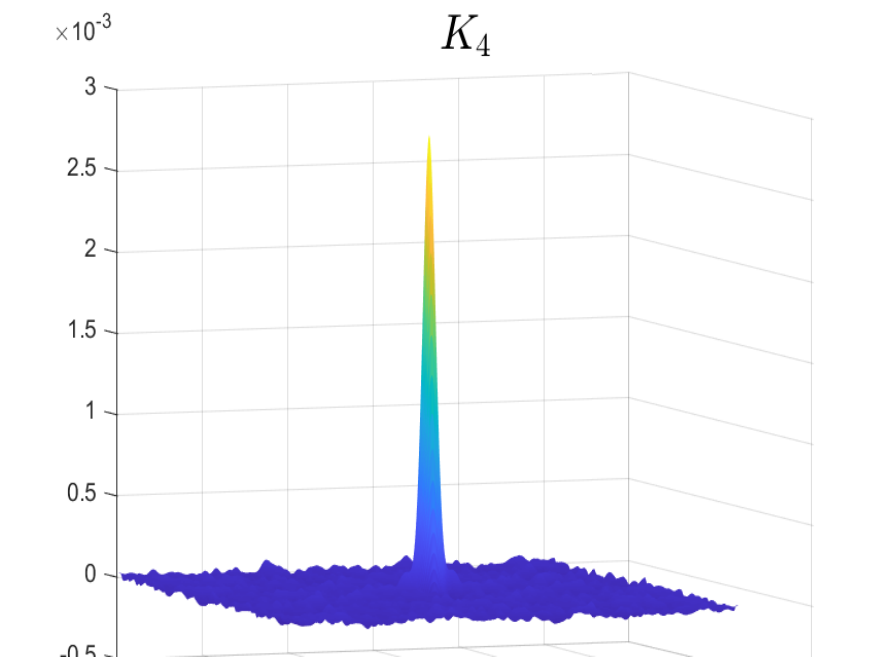}
\end{minipage}
    \begin{minipage}[c]{0.49\textwidth}
   \includegraphics[scale = 0.4]{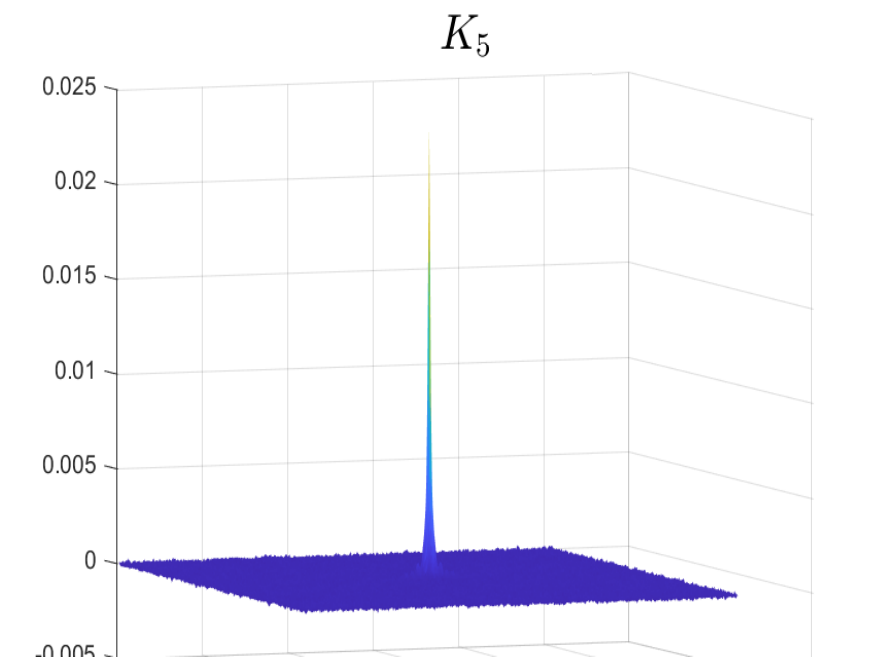}
\end{minipage}\\
 \begin{minipage}[c]{0.4\textwidth}
   \includegraphics[scale = 0.5]{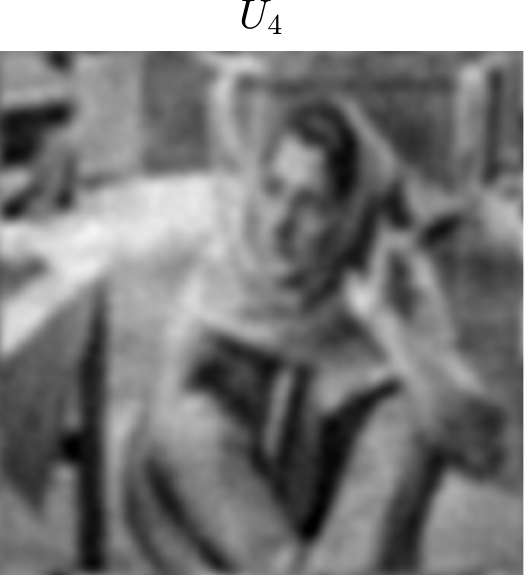}
\end{minipage}  \hspace{1.5cm}  \begin{minipage}[c]{0.4\textwidth}
   \includegraphics[scale = 0.5]{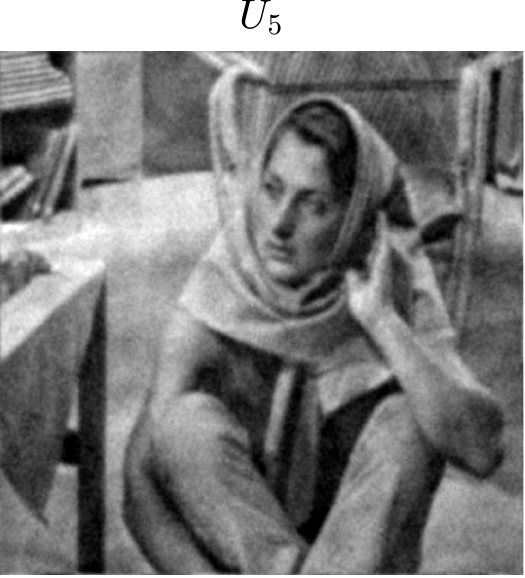}
\end{minipage}
    \caption{Left: reconstructed kernel and image at stopping index for data $f_1^\delta$. Right: reconstructed kernel and image at stopping index for data $f_2^\delta$.}
    \label{fig:blind_reconstructions}
\end{figure}

\begin{figure}
    \centering
\begin{minipage}[c]{0.45\textwidth}
    \includegraphics[scale = 0.5]{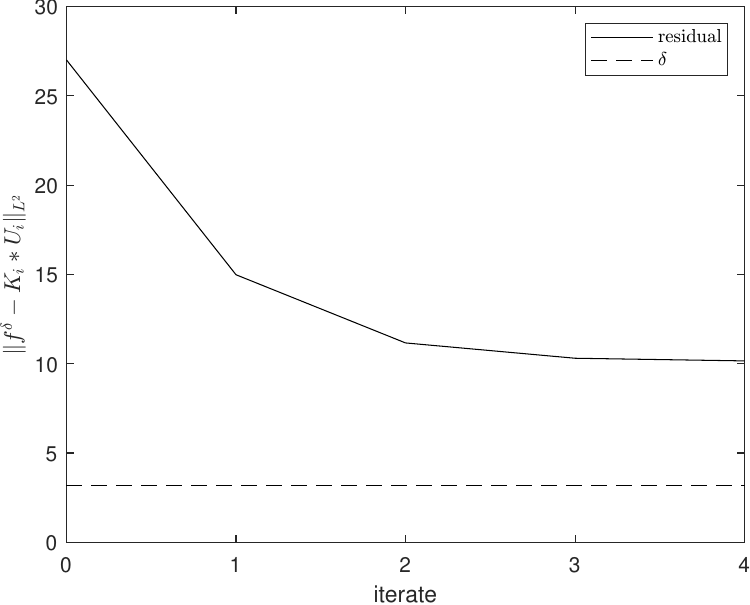}
\end{minipage}
\begin{minipage}[c]{0.45\textwidth}
    \includegraphics[scale = 0.5]{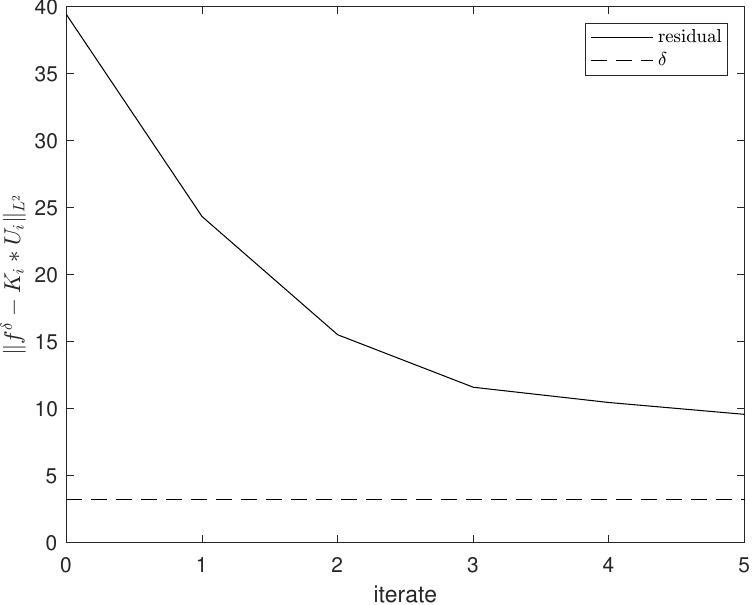}
\end{minipage}
    \caption{Left: Residual of the blind MHDM and noise level for $f_1^\delta$, Right: Residual of the blind MHDM and noise level for $f_2^\delta$.}
    \label{fig:H^r-H^s-blind_residual}
\end{figure}

\begin{figure}
    \centering
\begin{minipage}[c]{0.3\textwidth}
  \hspace{1cm} \includegraphics[scale = 0.4]{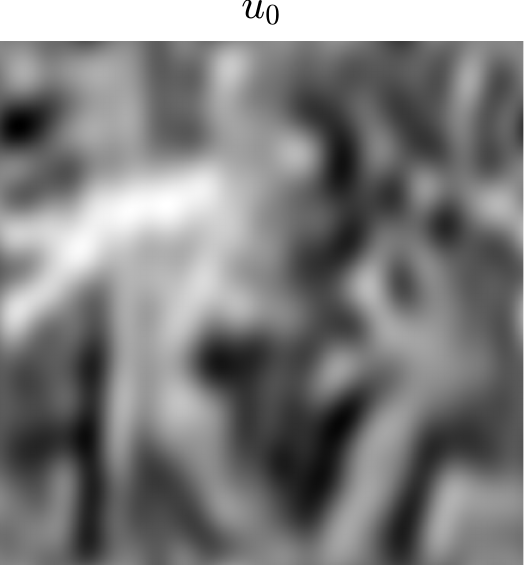}
\end{minipage}
\begin{minipage}[c]{0.3\textwidth}  
\hspace{0cm}\includegraphics[scale = 0.4]{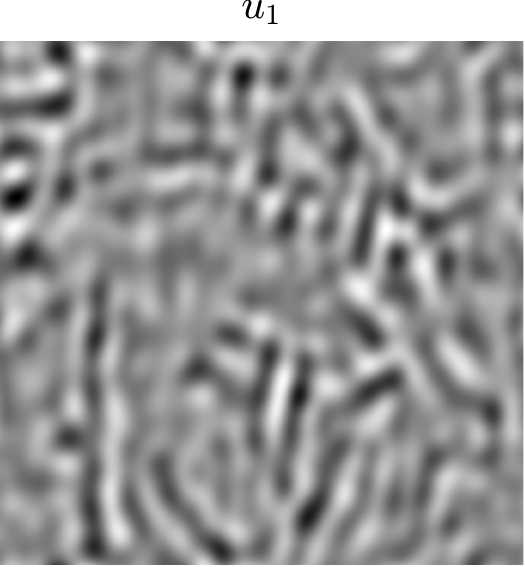}
\end{minipage}
\begin{minipage}[c]{0.3\textwidth}
  \hspace{-1.3cm} \includegraphics[scale = 0.4]{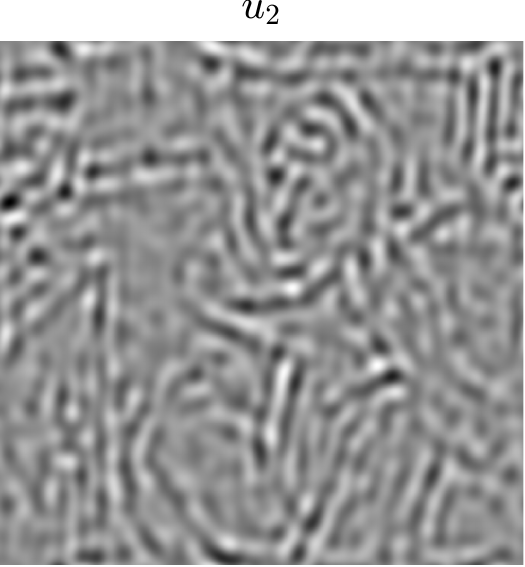}
\end{minipage}\\
\hspace{1.3cm}\begin{minipage}[c]{0.3\textwidth}
 \hspace{0.5cm}\includegraphics[scale = 0.4]{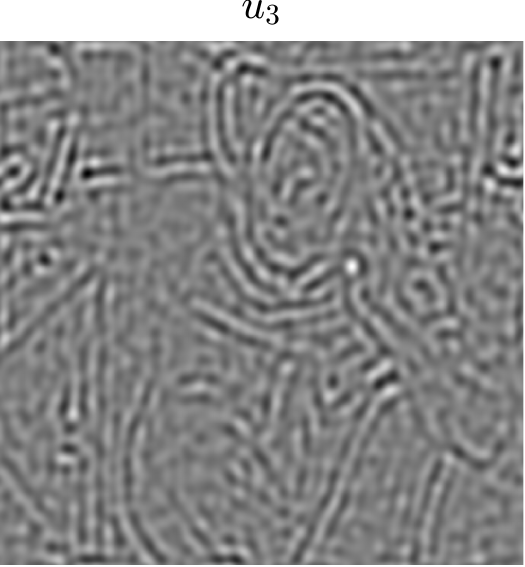}
\end{minipage}
\begin{minipage}[c]{0.3\textwidth}
  \hspace{-0.75cm} \includegraphics[scale = 0.4]{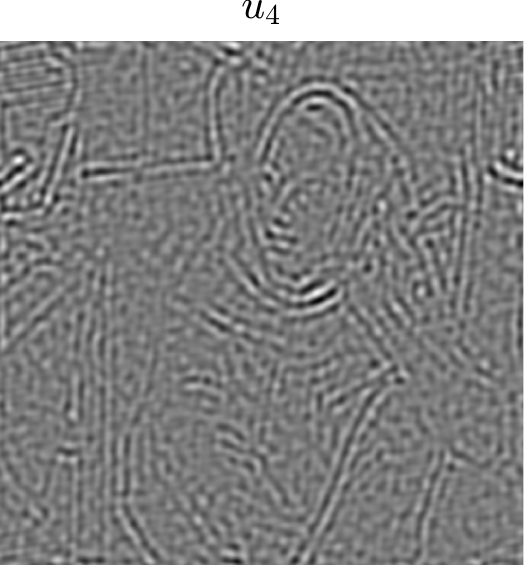}
\end{minipage}
\begin{minipage}[c]{0.3\textwidth}
\hspace{-1.9cm} \includegraphics[scale = 0.4]{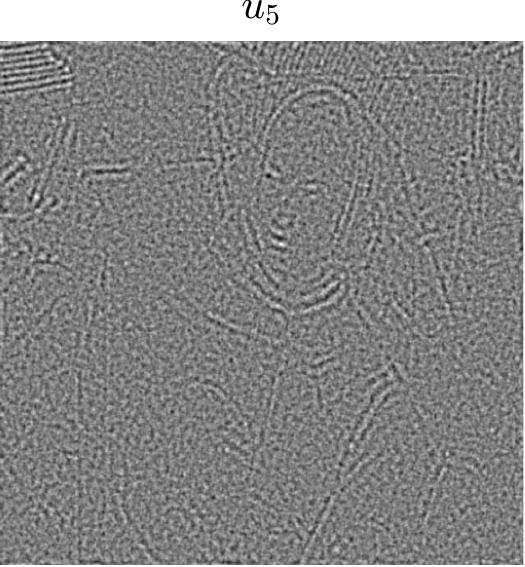}
\end{minipage}
 \caption{Scale decomposition $u_n$ obtained from $f_2^\delta$.}
    \label{fig:H^r-H^s-scales_image}
\end{figure}

\begin{figure}
    \centering
\begin{minipage}[c]{5cm}
\hspace{-1cm}\includegraphics[scale = 0.4]{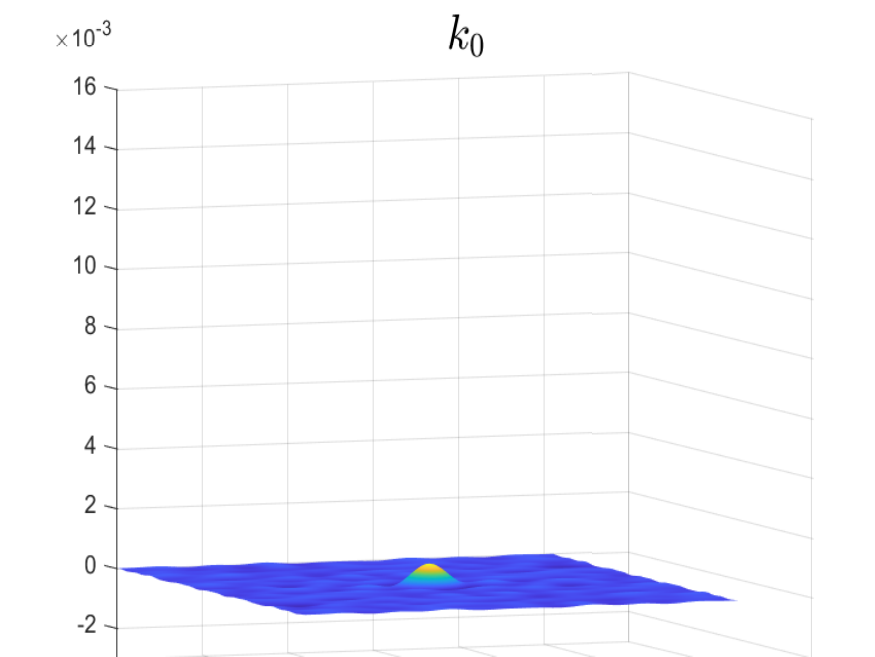}
\end{minipage}
\begin{minipage}[c]{5cm}
 \hspace{-1.2cm} \includegraphics[scale = 0.4]{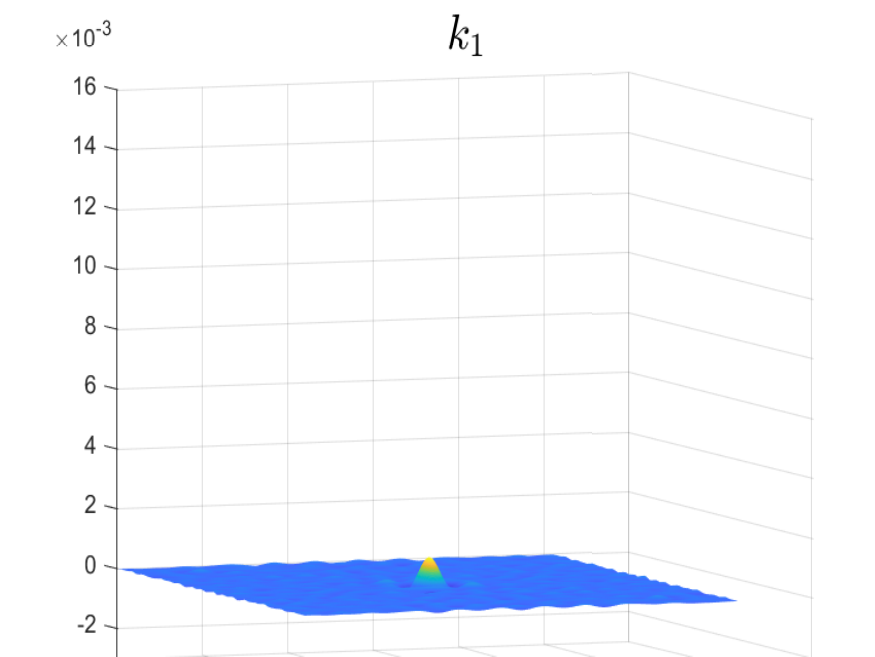}
\end{minipage}
\begin{minipage}[c]{5cm}
\hspace{-1.6cm} \includegraphics[scale = 0.4]{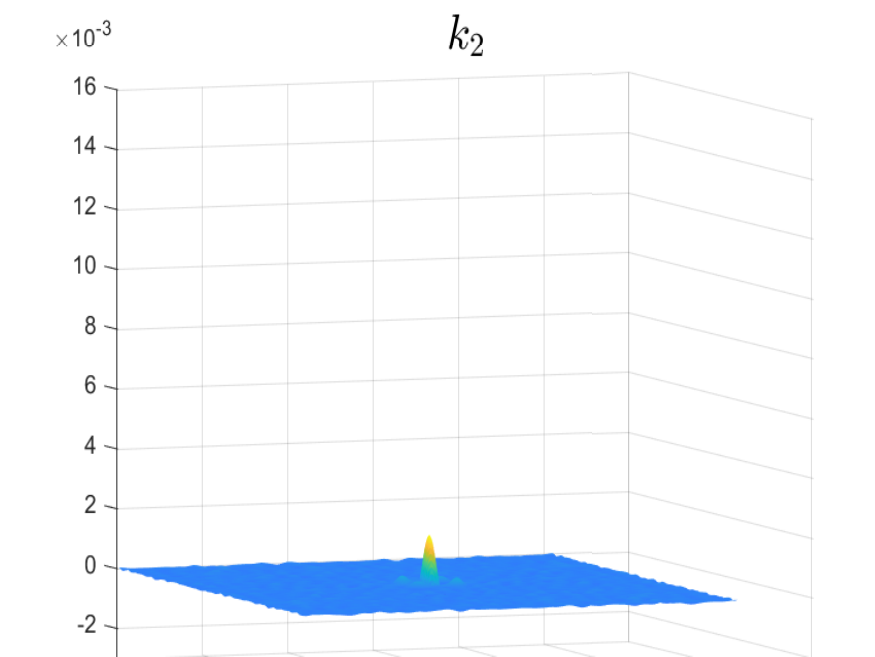}
\end{minipage}
\\
\begin{minipage}[c]{5cm}
\hspace{-1cm}\includegraphics[scale = 0.4]{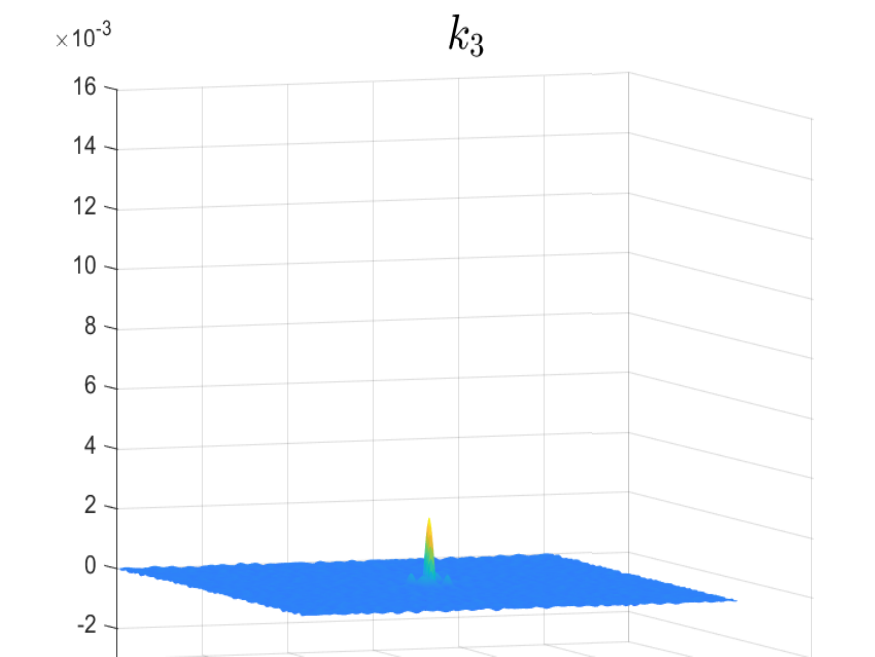}
\end{minipage}
\begin{minipage}[c]{5cm}
 \hspace{-1.2cm} \includegraphics[scale = 0.4]{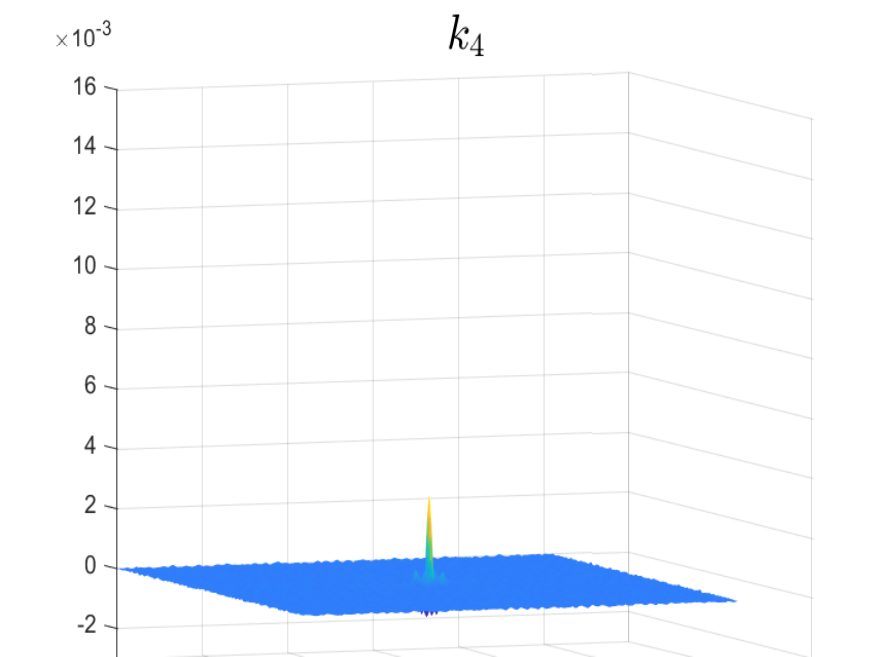}
\end{minipage}
\begin{minipage}[c]{5cm}
\hspace{-1.6cm} \includegraphics[scale = 0.4]{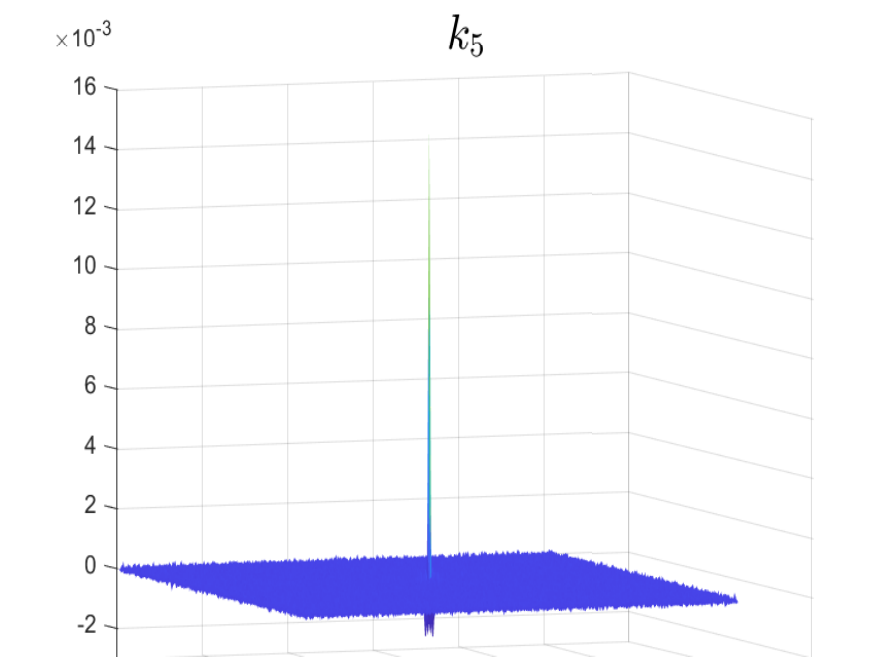}
\end{minipage}
 \caption{Scale decomposition $k_n$ obtained from $f_2^\delta$.}
    \label{fig:H^r-H^s-scales_kernel}
\end{figure}

\begin{figure}
    \centering
\begin{minipage}[c]{5cm}
\hspace{-1cm}\includegraphics[scale = 0.4]{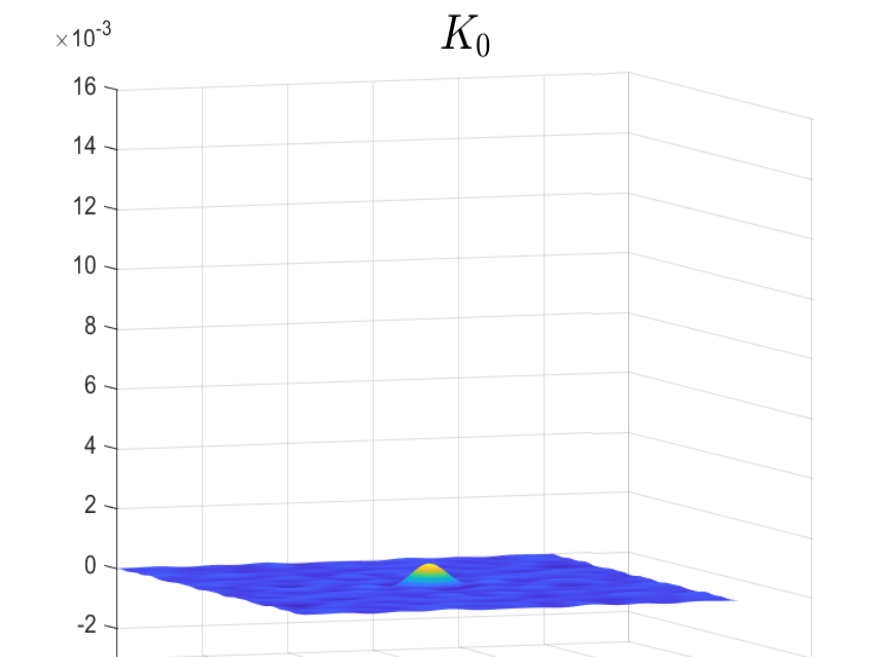}
\end{minipage}
\begin{minipage}[c]{5cm}
 \hspace{-1.2cm} \includegraphics[scale = 0.4]{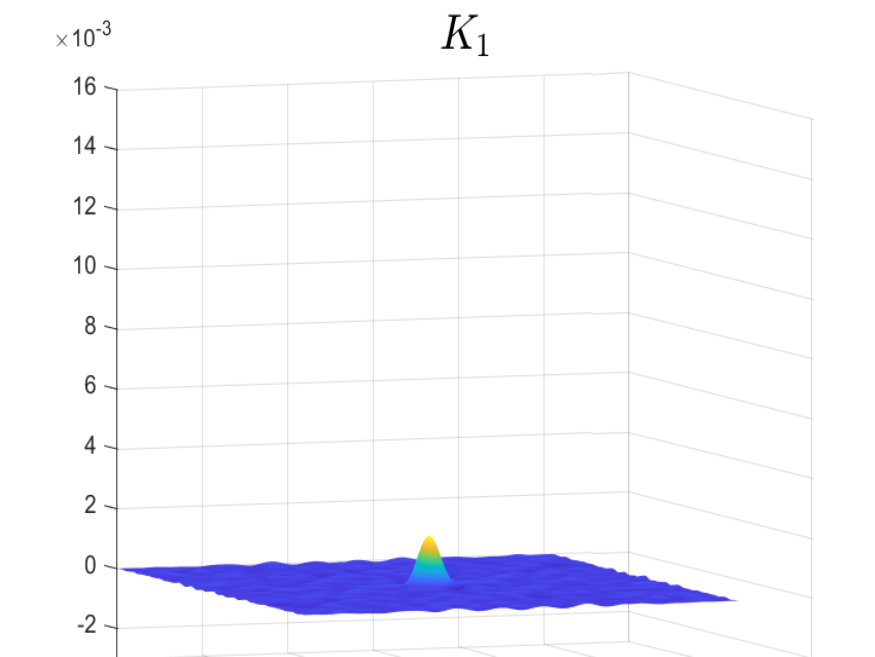}
\end{minipage}
\begin{minipage}[c]{5cm}
\hspace{-1.6cm} \includegraphics[scale = 0.4]{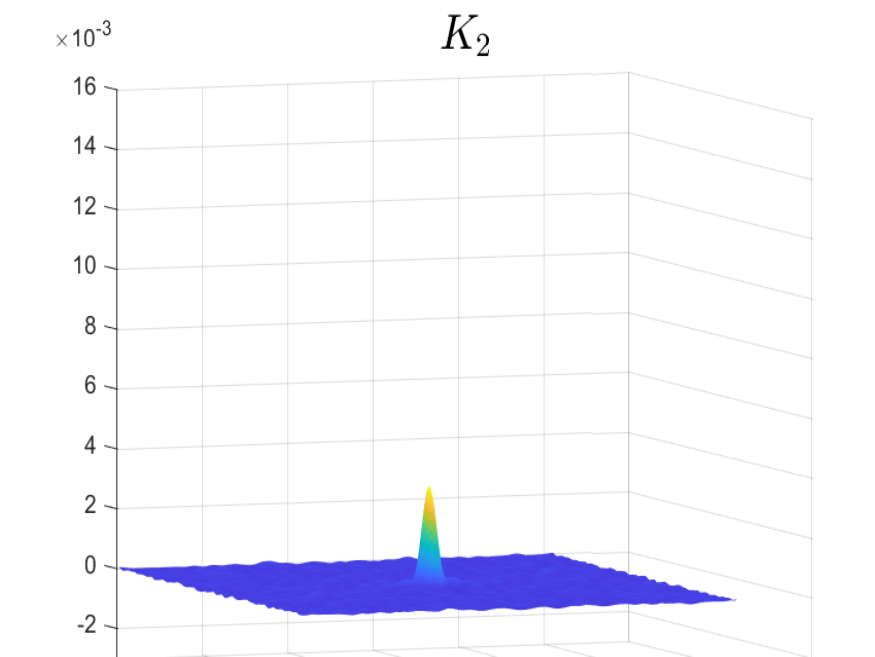}
\end{minipage}
\\
\begin{minipage}[c]{5cm}
\hspace{-1cm}\includegraphics[scale = 0.4]{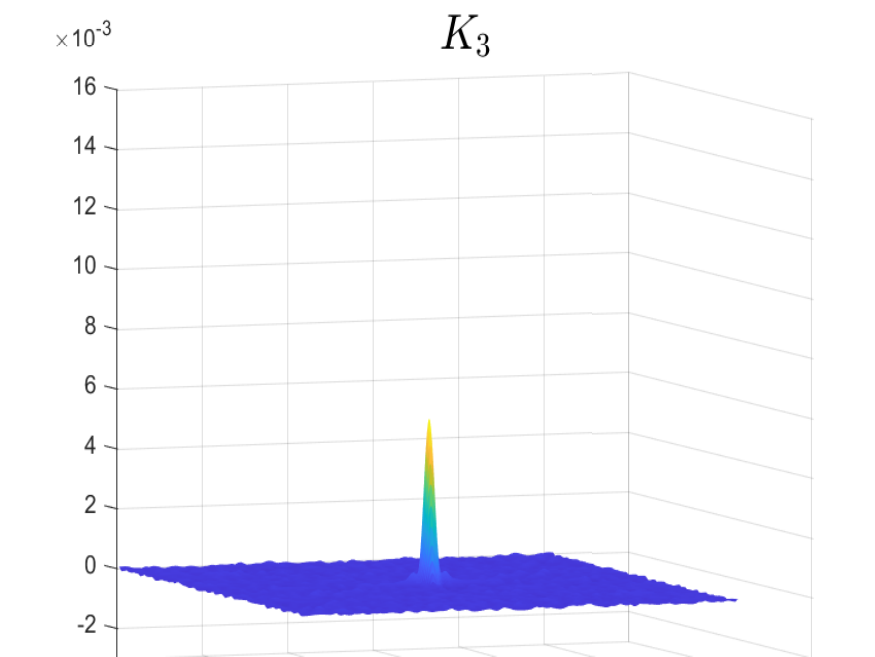}
\end{minipage}
\begin{minipage}[c]{5cm}
 \hspace{-1.2cm} \includegraphics[scale = 0.4]{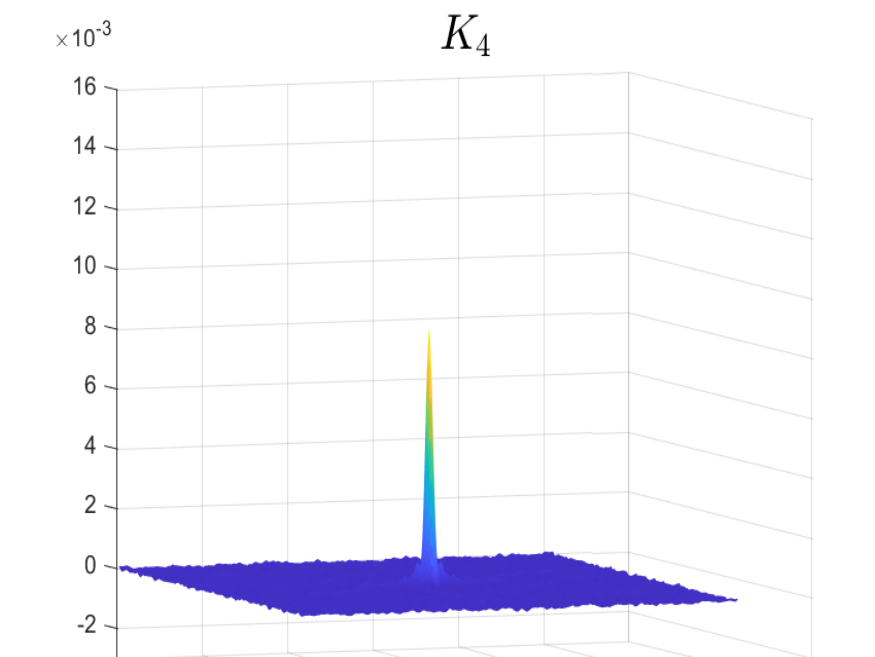}
\end{minipage}
\begin{minipage}[c]{5cm}
\hspace{-1.6cm} \includegraphics[scale = 0.4]{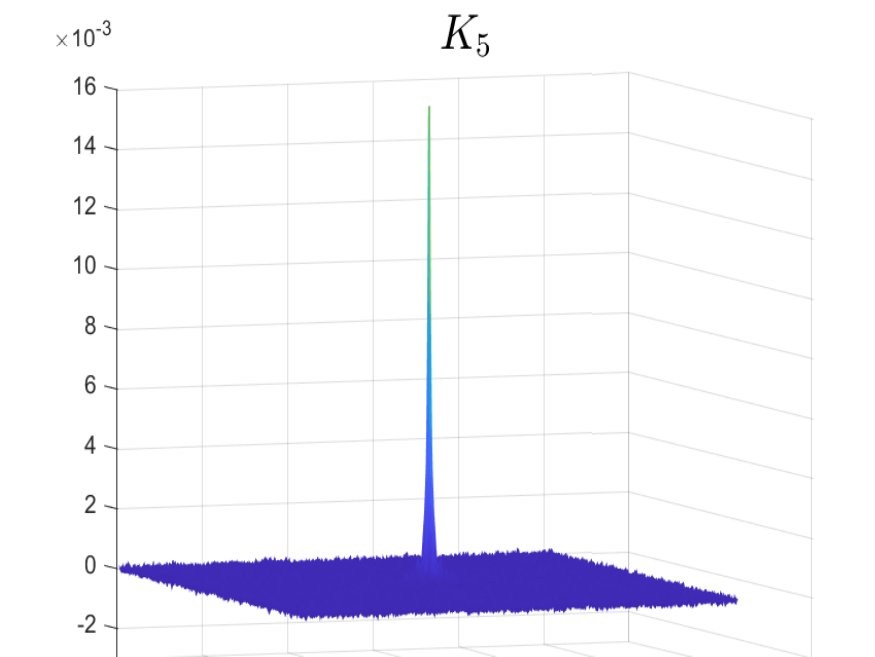}
\end{minipage}
 \caption{Iterates $K_n = \sum\limits_{i = 0}^n k_i$ obtained from $f_2^\delta$.}
    \label{fig:H^r-H^s-iterates_kernel}
\end{figure}

Moreover,  the experiments suggest that a suboptimal choice for the initial ratio leads to a worse kernel reconstruction while, upon an affine rescaling  of the grayscale values, the reconstructed images are visually still good. This is illustrated in Figure \ref{fig:images_different_parameters}, where the reconstructed kernels and images (with the pixel values rescaled) for different initial parameters $\lambda_0$ and $\mu_0$ are shown. We observe that visually all images are indistinguishable and approximate the true image well. Similarly, the corresponding kernels are structurally similar to the real one, but their numerical values are very disproportional.

\begin{figure}
    \centering
\begin{minipage}[c]{0.24\textwidth}
 \hspace{-1cm} \includegraphics[scale = 0.3]{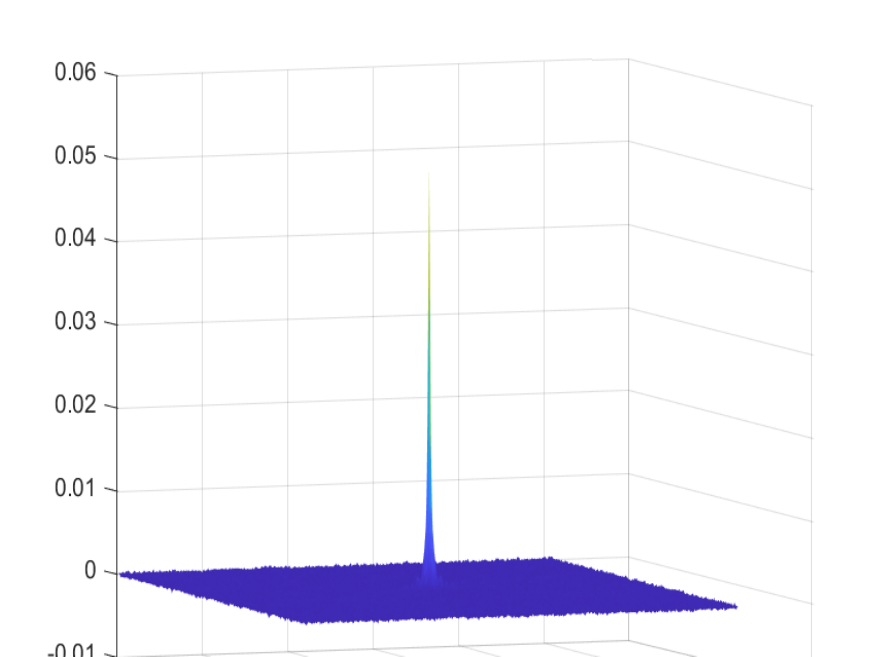}
\end{minipage}
\begin{minipage}[c]{0.24\textwidth}
\hspace{-1cm} \includegraphics[scale = 0.3]{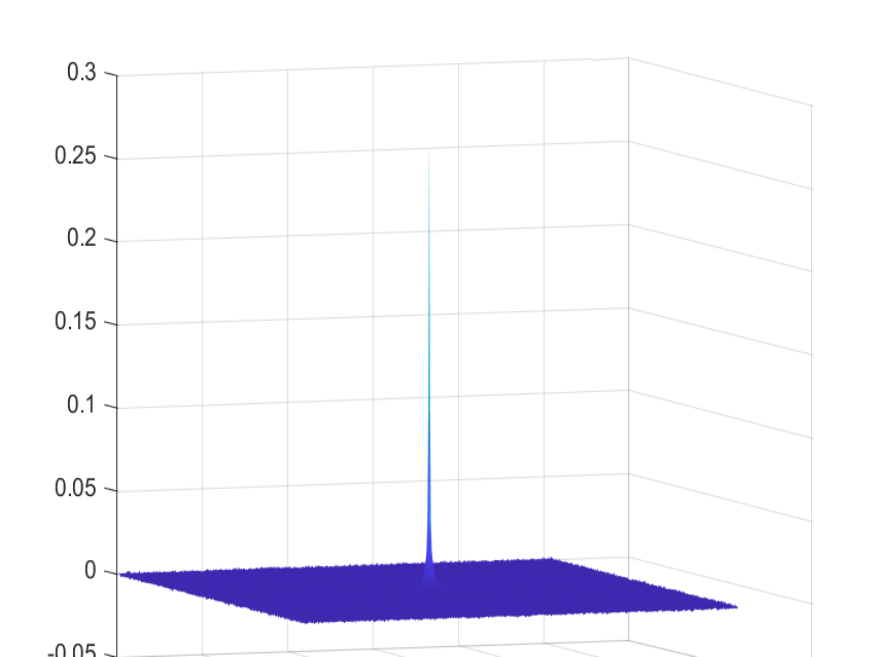}
\end{minipage}
\begin{minipage}[c]{0.24\textwidth}
\hspace{-1cm} \includegraphics[scale = 0.3]{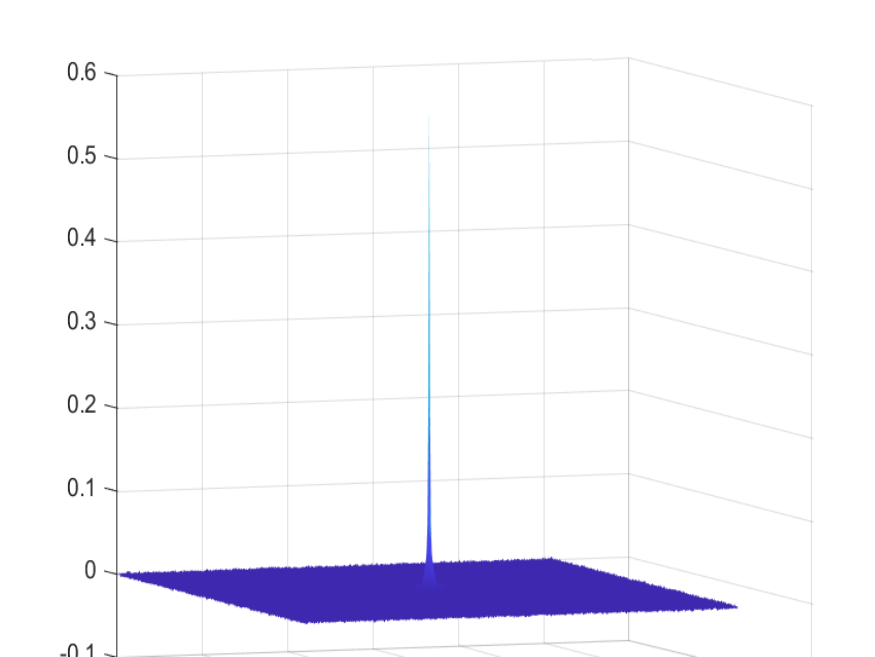}
\end{minipage}
\begin{minipage}[c]{0.24\textwidth}
 \hspace{-1cm}  \includegraphics[scale = 0.3]{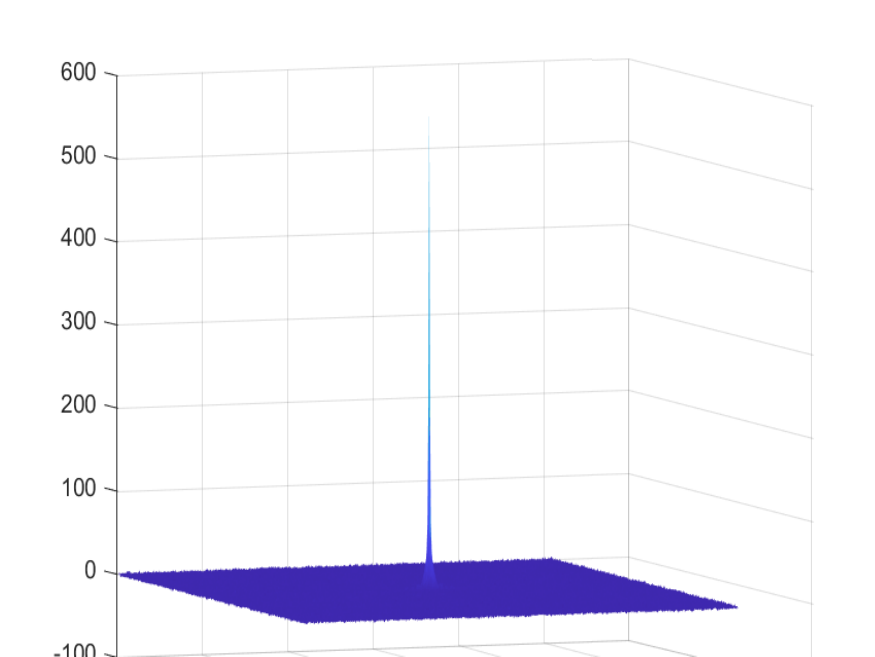}
\end{minipage}\\
    
\begin{minipage}[c]{0.24\textwidth}
  \hspace{-0.5cm}  \includegraphics[scale = 0.4]{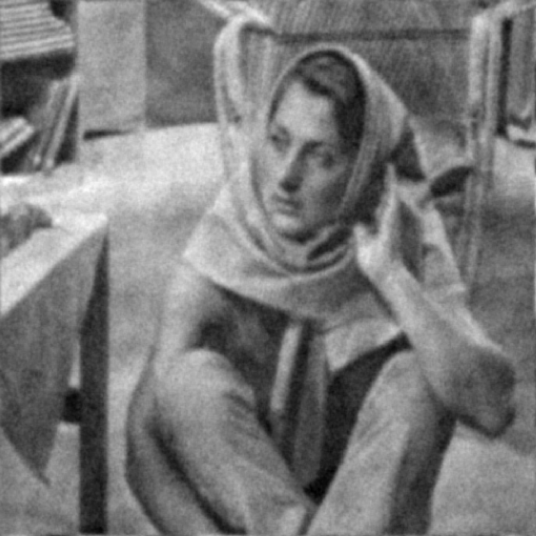}
\end{minipage}
\begin{minipage}[c]{0.24\textwidth}
  \hspace{-0.5cm}   \includegraphics[scale = 0.4]{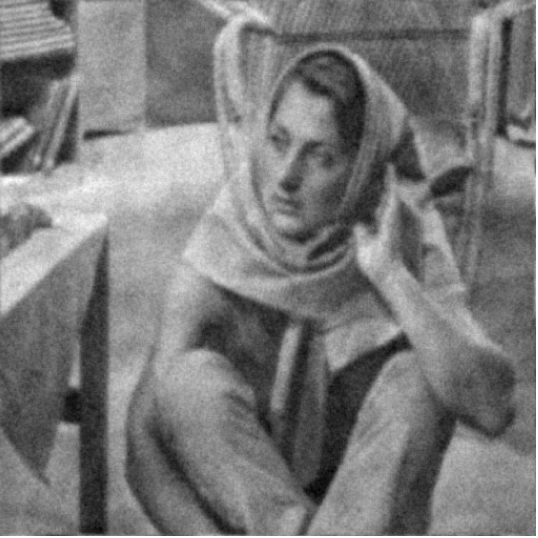}
\end{minipage}
\begin{minipage}[c]{0.24\textwidth}
 \hspace{-0.5cm}    \includegraphics[scale = 0.4]{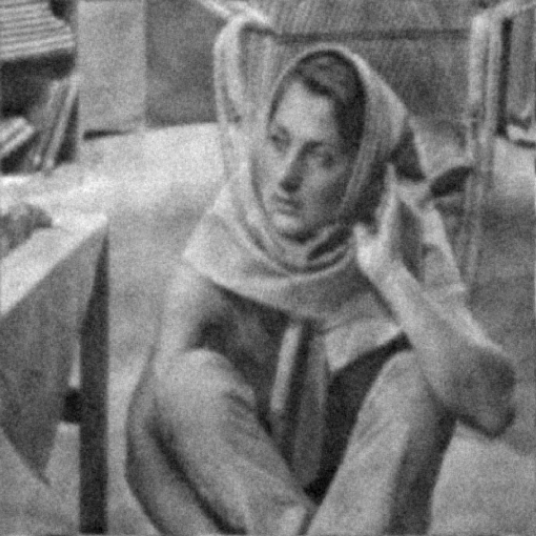}
\end{minipage}
\begin{minipage}[c]{0.24\textwidth}
 \hspace{-0.5cm}    \includegraphics[scale = 0.4]{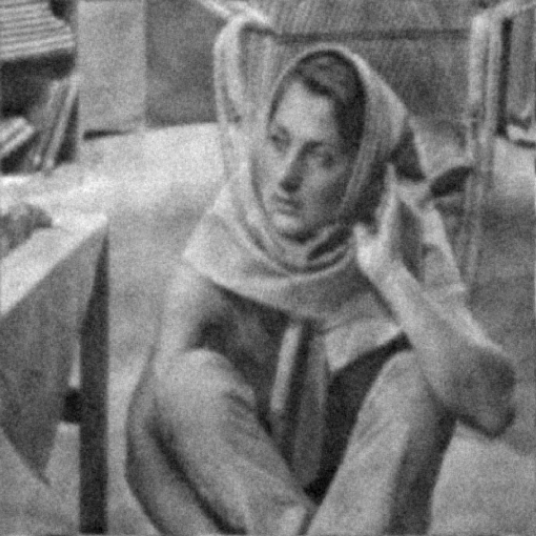}
\end{minipage}
    \caption{Kernel (top) and image (bottom) reconstructions obtained by running the blind-deconvolution MHDM with different initial parameters and rescaling the pixel values. From left to right: $\lambda_0 = 2\times 10^{-3}$ and $\mu_0 = 6.3\times 10^5$, $\lambda_0 = 1.4\times 10^{-4}$ and $\mu_0 = 1\times 10^3$,  $\lambda_0 = 2\times 10^{-3}$ and $\mu_0 = 1\times 10^3$,  $\lambda_0 = 2\times 10^{3}$ and $\mu_0 = 1\times 10^{-3}$. }
    \label{fig:images_different_parameters}
\end{figure}

\subsection{Experiment 1: Comparison  of blind MHDM vs. non-blind MHDM}

We compare the proposed blind MHDM to a non-blind version of the MHDM  (see for instance \cite{VeseDeblurring, MultiscaleRefinementImaging}) that uses the penalty term $J = \norm{\cdot}_{H^r}^2$ for the image. Instead of including the reconstruction of the blurring kernel in the method, we simply make a guess in the non-blind MHDM and stop the iteration once the discrepancy principle is satisfied. The non-blind MHDM therefore only requires the choice of one initial parameter $\lambda_0$, which we choose to be the same as for the blind deconvolution MHDM and also decrease according to the rule $\lambda_n = 4^{-n} \lambda_0$. To compare the two methods, we test the non-blind MHDM for $1000$ centered Gaussian kernels with variances ranging between $1$ and $12$ as guessed kernels. 

In Figures~\ref{fig:H^r-H^s-comp_non_blind_image} and \ref{fig:H^r-H^s-comp_non_blind_kernel}, we compare 
the performance of the blind MHDM with 
the non-blind MHDM algorithm in terms of 
 the error measures PSNR, SSIM 
(Figure~\ref{fig:H^r-H^s-comp_non_blind_image}), and the 
$L^2$-error of the kernel (Figure~\ref{fig:H^r-H^s-comp_non_blind_kernel}). 
We consider either the 
data $f_1^\delta$ (figures on the left-hand side) or
$f_2^\delta$ (figures on the right-hand side). In each of the plots, the values on the $x$-axis correspond to the guessed variance $\sigma$ of the Gaussian kernel used for the non-blind MHDM. The full line corresponds to the values of the resulting respective error measures, whereas the constant dashed line represents the value of the blind MHDM. Clearly, the quality of the non-blind algorithms depends on 
the correctness of the choice of $\sigma$. 
For the observation $f_1^\delta$, which corresponds to blurring with a single Gaussian kernel, the non-blind MHDM expectedly outperforms the blind MHDM only for kernel guesses that are similar to the true kernel. We observe similar behavior for the case $f_2^\delta$ with multiple Gaussian kernels used in the blurring. 
However, let us point out that, especially 
for the error measures PSNR and SSIM, 
the non-blind method has only a rather modest advantage and only in the case 
when the “guessed” $\sigma$ is close to 
the true one, while in case of a wrong 
guess, the non-blind method can go wrong 
quite dramatically as can be seen from the experiments with observation $f_2^\delta$ (the bottom line in Figure \ref{fig:H^r-H^s-comp_non_blind_image}). To further test this information, we repeat the experiment for $16$ different test images, where each of those is corrupted by $3$ different blurring kernels: Kernel $1$ is a single centered Gaussian kernel with variance  $8$ representing a strong blur, Kernel $2$ is the convex combination of $4$ different Gaussian kernels with variances ranging from $1$ to $5$ used in the previous experiments, while Kernel $3$ is a centered Gaussian kernel with variance $2$ and represents a mild blur. Additionally, to all blurred images we add Gaussian noise with respective variances $0.0004, 0.00004, 0.003$. Thus we consider a total of $144$ different test images. We then compute approximate solutions of the blind deblurring problem for these images using the MHDM with discrepancy principle for $\tau  =\sqrt{1.001}$. Thereafter, we again run $1000$ non-blind versions of the MHDM using $1000$ ``guessed'' Gaussian kernels with variances between $1$ and $12$. In Table \ref{tab:avg_noise_non_blind}  we report the average ratio of the PSNR and SSIM values between the image obtained via blind MHDM and the best recontructed image froom the non-blind method with respect to the noise or blurring kernel. 

\begin{table}[H]
\centering
\begin{tabular}{|c|c|c|}
     \hline $\sigma_{noise}$  & $\frac{\text{PSNR}_\text{blind}}{\text{PSNR}_{\text{guess}}}$ & $\frac{\text{SSIM}_\text{blind}}{\text{SSIM}_{\text{guess}}}$ \\ \hline
     0.0004&0.870&0.903\\4e-05&0.855&0.892\\0.003&0.892&0.906\\\hline
     \end{tabular}
     \begin{tabular}{|c|c|c|}
     \hline kernel  & $\frac{\text{PSNR}_\text{blind}}{\text{PSNR}_{\text{guess}}}$& $\frac{\text{SSIM}_\text{blind}}{\text{SSIM}_{\text{guess}}}$ \\ \hline
     1&0.907&0.931\\2&0.875&0.905\\3&0.835&0.865\\\hline
     \end{tabular}
    \caption{Average ratio of PSNR and SSIM between the images obtained via blind and non-blind deconvolution MHDM sorted by noise level and blurring kernel respcetively.}
    \label{tab:avg_noise_non_blind}
\end{table}


We note that, on average, the blind MHDM performs at most $13\%$ worse than the optimal non-blind MHDM. These results indicate that the blind MHDM is a  robust method 
that produces reasonable approximations.
In particular, if the  blurring kernel 
cannot be estimated with high accuracy, 
then the blind MHDM is superior to non-blind approaches.

\begin{figure}
    \centering
\begin{minipage}[c]{0.45\textwidth}
    \includegraphics[scale = 0.5]{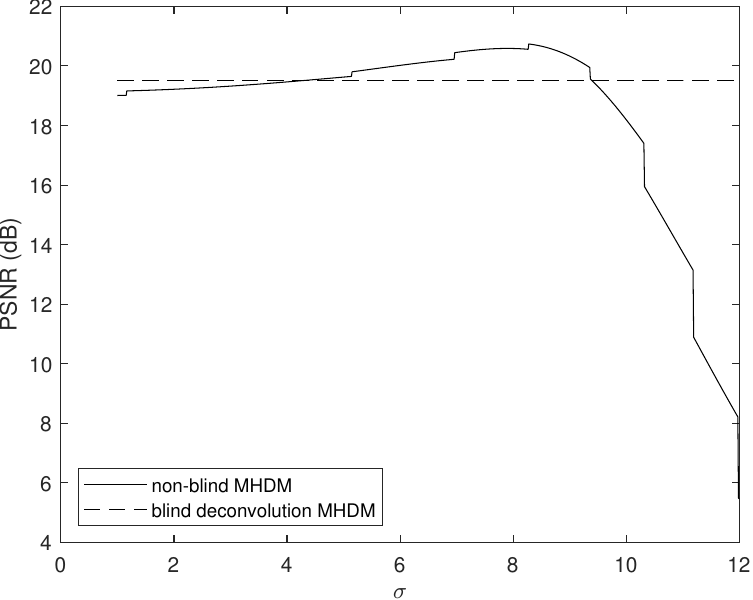}
\end{minipage}
\begin{minipage}[c]{0.45\textwidth}
    \includegraphics[scale = 0.5]{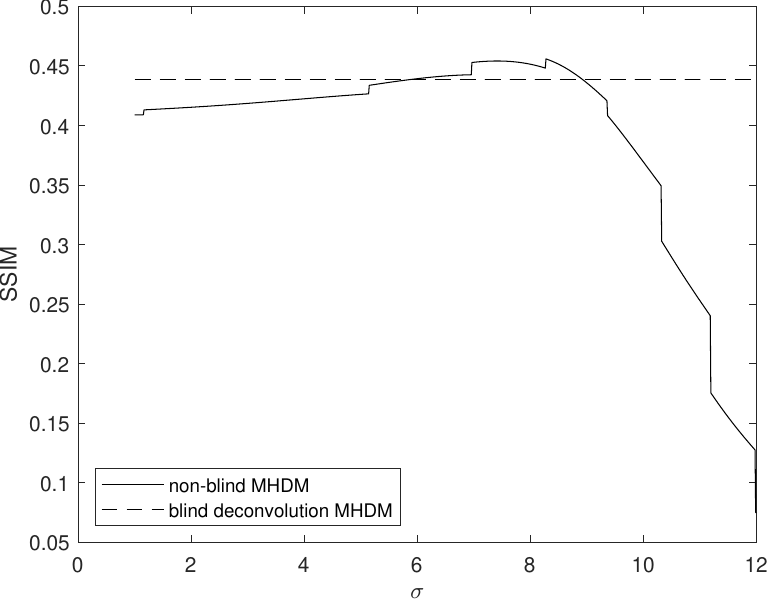}
\end{minipage}\\
\begin{minipage}[c]{0.45\textwidth}
    \includegraphics[scale = 0.5]{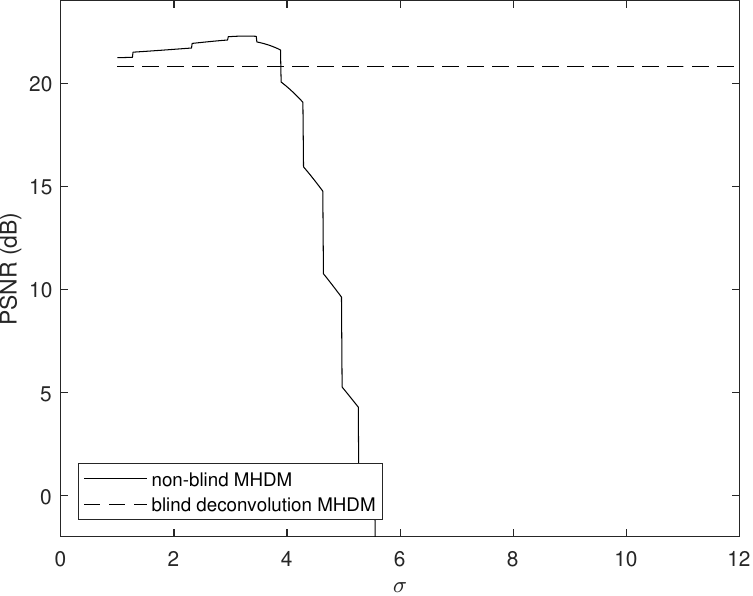}
\end{minipage}
\begin{minipage}[c]{0.45\textwidth}
    \includegraphics[scale = 0.5]{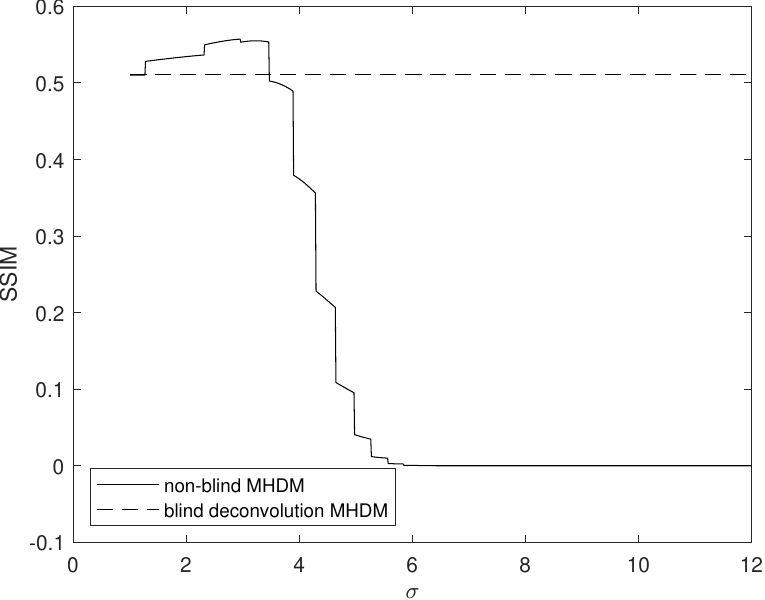}
\end{minipage}
    \caption{Top left: PSNR values for non-blind MHDM with guessed kernel for different guesses of the kernel and data $f_1^\delta$,\\ Top right: SSIM values for non-blind MHDM with guessed kernel for different guesses of the kernel and data $f_1^\delta$,\\ Bottom left: PSNR values for non-blind MHDM with guessed kernel for different guesses of the kernel and data $f_2^\delta$,\\ Bottom right: SSIM values for non-blind MHDM with guessed kernel for different guesses of the kernel and data $f_2^\delta$.}
    \label{fig:H^r-H^s-comp_non_blind_image}
\end{figure}

\begin{figure}
    \centering
\begin{minipage}[c]{0.45\textwidth}
    \includegraphics[scale = 0.5]{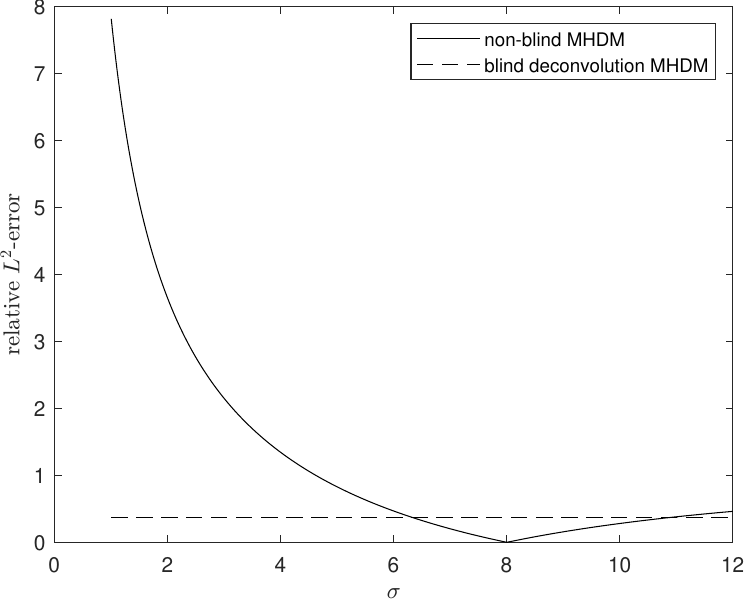}
\end{minipage}
\begin{minipage}[c]{0.45\textwidth}
    \includegraphics[scale = 0.5]{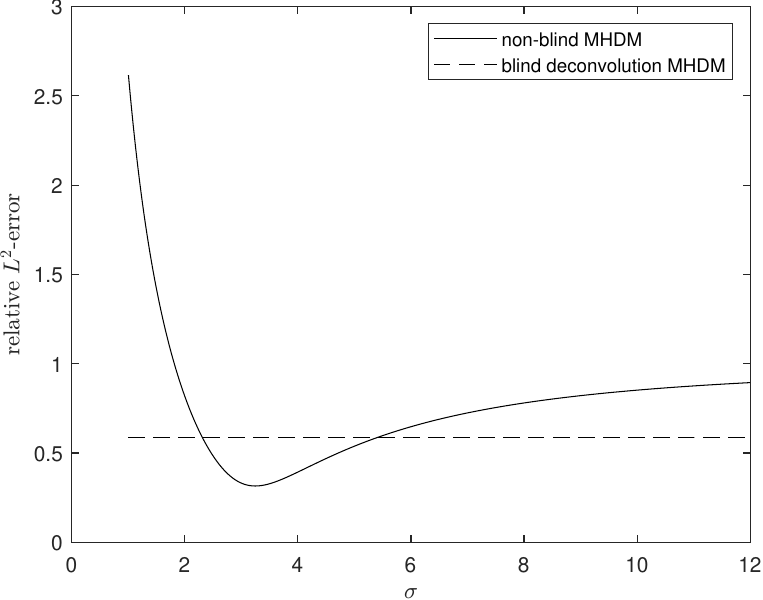}
\end{minipage}
    \caption{Comparison of the relative $L^2$ error of the guessed kernel with the relative error of the kernel reconstructed by the blind MHDM with observation $f_1^\delta$ (left) and $f_2^\delta$ (right).}
    \label{fig:H^r-H^s-comp_non_blind_kernel}
\end{figure}

\subsection{Experiment 2: Comparison blind MHDM vs variational blind deconvolution}

We compare the blind deconvolution MHDM with a single-step variational regularization as in \eqref{eq:H^r-H^k_initial}. To achieve a fair comparsion of the methods, we optimize the regularization parameters. For the single-step regularization we proceed as follows: For a given ratio of regularization parameters $\frac{\mu_0}{\lambda_0}$, we compute regularized solutions of the blind-deconvolution problem by a grid search. As for the MHDM, we start this grid search by computing minimizers for an initial choice of parameters $(\lambda,\mu)$. The consecutive iterates $(u_n^{var},k_n^{var})$ are then obtained by computing solutions of \eqref{eq:H^r-H^k_initial} with parameters $(\lambda_n,\mu_n) = (4^n \lambda_0,4^n\mu_0)$. This procedure is terminated once the discrepancy principle \eqref{eq:discrepancy_principle} with $\tau = \sqrt{1.001}$ is satisfied. Therefore, the grid search and the MHDM are stopped according to the same rule which allows to investigate the effect of the multiscale decomposition. We then optimize both methods with respect to the ratio $\frac{\mu_0}{\lambda_0}$ of the initial parameters. Here we consider optimal those parameters which maximize the PSNR value of the reconstructed image at the iterate satisfying the stopping rule. In order to observe the effects of the multiscale decomposition, we additionally require that both methods run for at least 10 iterations until the discrepancy principle is met. This way, we ensure that there are enough iterations such that the multiscale effects    are observable on the reconstructions.
For the experiments, we consider the same $144$ test images as in the previous section.
The average PSNR and SSIM values of the reconstructed images and the average relative $L^2$-error of the reconstructed kernels, as well as the corresponding standard deviations for both methods, are shown in Table \ref{tab:avg_noise} (averages and standard deviations for the different levels of noise) and Table \ref{tab:avg_blur} (averages and standard deviations for the different blurring kernels). We observe that, on average, the MHDM performs slightly better with respect to all considered metrics, apart from the average PSNR for noise with variance $0.0004$. Based on a closer inspection of the experimental data, we notice that the advantage of the MHDM is, in part,  due to the MHDM terminating more frequently with reconstructions whose residual is closer to the stopping threshold $\tau\delta^2$. This is illustrated in Tables \ref{tab:avg_noise_cut} and \ref{tab:avg_blur_cut}, where we omit those experiments in which the residuals of the MHDM and the variational grid search at the stopping index differ by more than $0.1\tau \delta^2$. We notice that, by neglecting  these ``outliers'', both methods perform even more similarly with respect to the quality indicators. By our numerical tests we thus conclude that  the MHDM produces qualitatively similar approximations of the true image and kernel, and typically terminates with a residual closer to the stopping threshold. This means the MHDM is more robust than a grid search for single-step blind-deblurring.\\ 
Furthermore, we observe for both methods that the optimal ratio $\frac{\mu_0}{\lambda_0}$ depends more on the true image and less on the kernel used to obtain the blurred image. Additionally, for a given image, there is only a slight difference between the optimal ratios of initial parameters   for the MHDM and the grid search in the single-step approach.

\begin{table}[H]
    \centering
\resizebox{\textwidth}{!}{\begin{tabular}{|c|cc|cc|cc|}
     \hline $\sigma_{noise}$  &$\text{PSNR}_{\text{MHDM}}$ &$\text{PSNR}_{\text{var}}$ & $\text{SSIM}_{\text{MHDM}}$ & $\text{SSIM}_{\text{var}}$ &$\text{err}_{\text{MHDM}}$ &$\text{err}_{\text{var}}$ \\ \hline
     0.0004&20.843 (2.524)&\textbf{20.843} (2.527)&\textbf{0.594} (0.124)&0.593 (0.125)&\textbf{0.492} (0.135)&0.496 (0.133)\\0.00004&\textbf{21.128} (2.641)&21.122 (2.639)&\textbf{0.624} (0.131)&0.623 (0.131)&\textbf{0.506} (0.180)&0.518 (0.219)\\0.003&\textbf{20.444} (2.344)&20.393 (2.352)&\textbf{0.554} (0.123)&0.534 (0.128)&\textbf{0.590} (0.214)&0.703 (0.244)\\\hline
     \end{tabular}}
    \caption{Comparison of the MHDM and the single-step variational blind deblurring for the different noise variances with respect to the average PSNR and SSIM values for the reconstructed image and average relative $L^2$-error of the reconstructed kernel. The numbers in brackets are the corresponding standard deviations.}
    \label{tab:avg_noise}
\end{table}

\begin{table}[H]
    \centering
\resizebox{\textwidth}{!}{\begin{tabular}{|c|cc|cc|cc|}
     \hline kernel  &$\text{PSNR}_{\text{MHDM}}$ &$\text{PSNR}_{\text{var}}$ & $\text{SSIM}_{\text{MHDM}}$ & $\text{SSIM}_{\text{var}}$ &$\text{err}_{\text{MHDM}}$ &$\text{err}_{\text{var}}$ \\ \hline
     1&\textbf{19.280} (1.825)&19.273 (1.823)&\textbf{0.517} (0.129)&0.510 (0.127)&\textbf{0.473} (0.192)&0.535 (0.206)\\2&\textbf{21.359} (2.362)&21.354 (2.355)&\textbf{0.612} (0.112)&0.612 (0.113)&\textbf{0.569} (0.145)&0.582 (0.186)\\3&\textbf{21.775} (2.554)&21.732 (2.590)&\textbf{0.642} (0.110)&0.628 (0.128)&\textbf{0.547} (0.197)&0.602 (0.269)\\\hline
     \end{tabular}}
    \caption{Comparison of the MHDM and the single-step variational blind deblurring for the different blurring kernels with respect to the average PSNR and SSIM values for the reconstructed image and average relative $L^2$-error of the reconstructed kernel. The numbers in brackets are the corresponding standard deviations.}
    \label{tab:avg_blur}
\end{table}

\begin{table}[H]
    \centering
\resizebox{\textwidth}{!}{\begin{tabular}{|c|cc|cc|cc|}
     \hline $\sigma_{noise}$  &$\text{PSNR}_{\text{MHDM}}$ &$\text{PSNR}_{\text{var}}$ & $\text{SSIM}_{\text{MHDM}}$ & $\text{SSIM}_{\text{var}}$ &$\text{err}_{\text{MHDM}}$ &$\text{err}_{\text{var}}$ \\ \hline
     0.0004&20.843 (2.524)&\textbf{20.843} (2.527)&\textbf{0.594} (0.124)&0.593 (0.125)&\textbf{0.492} (0.135)&0.496 (0.133)\\0.00004&\textbf{21.144} (2.667)&21.135 (2.666)&\textbf{0.627} (0.130)&0.627 (0.131)&0.504 (0.181)&\textbf{0.500} (0.178)\\0.003&20.428 (2.437)&\textbf{20.437} (2.444)&\textbf{0.552} (0.127)&0.546 (0.127)&\textbf{0.594} (0.223)&0.657 (0.198)\\\hline
     \end{tabular}}
         \caption{Comparison of the MHDM and the single-step variational blind deblurring for the different noise variances with respect to the average PSNR and SSIM values for the reconstructed image and average relative $L^2$-error of the reconstructed kernel after removing ``outliers''. The numbers in brackets are the corresponding standard deviations.}
    \label{tab:avg_noise_cut}
\end{table}

\begin{table}[H]
    \centering
\resizebox{\textwidth}{!}{\begin{tabular}{|c|cc|cc|cc|}
     \hline kernel  &$\text{PSNR}_{\text{MHDM}}$ &$\text{PSNR}_{\text{var}}$ & $\text{SSIM}_{\text{MHDM}}$ & $\text{SSIM}_{\text{var}}$ &$\text{err}_{\text{MHDM}}$ &$\text{err}_{\text{var}}$ \\ \hline
     1&\textbf{19.280} (1.825)&19.273 (1.823)&\textbf{0.517} (0.129)&0.510 (0.127)&\textbf{0.473} (0.192)&0.535 (0.206)\\2&\textbf{21.380} (2.383)&21.372 (2.377)&\textbf{0.615} (0.111)&0.615 (0.112)&0.568 (0.146)&\textbf{0.564} (0.143)\\3&21.880 (2.633)&\textbf{21.897} (2.636)&0.648 (0.112)&\textbf{0.649} (0.113)&0.547 (0.206)&\textbf{0.547} (0.204)\\\hline
     \end{tabular}}
         \caption{Comparison of the MHDM and the single-step variational blind deblurring for the different blurring kernels with respect to the average PSNR and SSIM values for the reconstructed image and average relative $L^2$-error of the reconstructed kernel after removing ``outliers''. The numbers in brackets are the corresponding standard deviations.}
    \label{tab:avg_blur_cut}
\end{table}

\section{Conclusion}
We introduce the Multiscale Hierarchical Decomposition Method for the blind deconvolution problem and show convergence of the residual in the noise-free case and then in the noisy data case by employing a discrepancy principle. To demonstrate the efficiency and behavior of the proposed method, we focus on employing fractional Sobolev norms as regularizers and develop a way to compute the appearing minimizers explicitly in a pointwise manner. We want to stress that in our experience, any variational approach to blind deconvolution should incorporate prior information on the expected blurring kernel. In our setting, this was done by enforcing a positivity constraint on the Fourier transform of the kernels, thus favoring, e.g.,  Gaussian structures. Numerical comparisons with a single-step variational method and a non-blind MHDM show that our approach produces comparable results,  in a more stable manner. Additionally, the scale decomposition of both  reconstructed kernel and image provides a meaningful interpretation of the involved iteration steps. For future work, this opens up the possibility to modify the method based on prior information of the underlying true solution. By using multiple penalty terms throughout the iteration, one could construct approximate solutions that admit different structures at different levels of detail. Nonetheless, we believe that at first a better understanding of iterates' convergence behavior is necessary to systematically refine the method. In future research, we aim to adapt the proposed method for parameter identification problems with unknown forward operator and the classification of blurring operators occuring in real applications.

\section{Acknowledgements}
We thank Michael Quellmalz (Technische Universität Berlin) for his valuable remarks and literature suggestion on the positivity of Fourier transforms. We want to thank the referee for the useful and constructive comments that greatly helped to improve the manuscript. This research was funded in part by the Austrian Science Fund (FWF) [10.55776/DOC78]. For open access purposes, the authors have applied a CC BY public copyright license to any author-accepted manuscript version arising from this submission.

\bibliographystyle{siamplain}
\bibliography{bibliography}
\end{document}